\documentclass[11pt]{article}
\usepackage{xcolor}	
	\newcommand{\blind}{0}
	
    \makeatletter
    \renewcommand\section{\@startsection {section}{1}{\z@}%
                                       {-3.5ex \@plus -1ex \@minus -.2ex}%
                                       {2.3ex \@plus.2ex}%
                                       {\normalfont\fontfamily{phv}\fontsize{16}{19}\bfseries}}
    \renewcommand\subsection{\@startsection{subsection}{2}{\z@}%
                                         {-3.25ex\@plus -1ex \@minus -.2ex}%
                                         {1.5ex \@plus .2ex}%
                                         {\normalfont\fontfamily{phv}\fontsize{14}{17}\bfseries}}
    \renewcommand\subsubsection{\@startsection{subsubsection}{3}{\z@}%
                                        {-3.25ex\@plus -1ex \@minus -.2ex}%
                                         {1.5ex \@plus .2ex}%
                                         {\normalfont\normalsize\fontfamily{phv}\fontsize{14}{17}\selectfont}}
    \makeatother
	\usepackage{amsmath}
	\usepackage{algorithm}
    \usepackage[noend]{algpseudocode}
    \usepackage{comment}
    \usepackage{bm}
    \usepackage{multirow}
    \usepackage{caption}
\usepackage{subcaption}
	\usepackage{graphicx}
	\usepackage{enumerate}
	\usepackage{adjustbox}
	\usepackage{natbib} %comment out if you do not have the package
	\usepackage{url} % not crucial - just used below for the URL
	    \usepackage{amsfonts, amsthm, latexsym, amssymb}
	    \usepackage{multirow}
	    \usepackage{array}
	    \newcolumntype{P}[1]{>{\centering\arraybackslash}p{#1}}
\begin{document}
		
			%%%%%%%%%%%%%%%%%%%%%%%%%%%%%%%%%%%%%%%%%%%%%%%%%%%%%%%%%%%%%%%%%%%%%%%%%%%%%%
		\def\spacingset#1{\renewcommand{\baselinestretch}%
			{#1}\small\normalsize} \spacingset{1}
		%%%%%%%%%%%%%%%%%%%%%%%%%%%%%%%%%%%%%%%%%%%%%%%%%%%%%%%%%%%%%%%%%%%%%%%%%%%%%%
		
		\if0\blind
		{
			\title{\bf Designing Drone Delivery Networks for Vaccine Supply Chain: A Case Study of Niger}
			\author{Maximilian Kolter $^a$, Sandra D. Eksigolu$^a$, Sarah Nurre Pinkley$^a$, Ruben A. Proano $^b$ \\
			$^a$ Industrial Engineering, University of Arkansas, Fayetteville, USA\\
             $^b$ Industrial and Systems Engineering, Rochester Institute of Technology, Rochester, USA}
			\date{}
			\maketitle
		} \fi
		
		\if1\blind
		{
            \title{\bf Designing Drone Delivery Networks for Vaccine Supply Chain: A Case Study of Niger}
			\author{Author information is purposely removed for double-blind review}
			\maketitle
		} \fi
		\bigskip
		
	\begin{abstract}
The focus of this research is to evaluate the use of drones for the delivery of pediatric vaccines in remote areas of low-income and low-middle-income countries. Delivering vaccines in these regions is challenging because of the inadequate road networks, and long transportation distances that make it difficult to maintain the cold chain's integrity during transportation. We propose a mixed-integer linear program to determine the location of drone hubs to facilitate the delivery of vaccines. The model considers the operational attributes of drones, vaccine wastage in the supply chain, cold storage and transportation capacities. We develop a case study using data from Niger to determine the impact of drone deliveries in improving vaccine availability in Niger. Our numerical analysis show a 0.71\% to 2.12\% increase in vaccine availability. These improvements depend on the available budget to build drone hubs and purchase drones, and the population density in the region of study.      
	\end{abstract}
			
	\noindent%
	\spacingset{1.5} % DON'T change the spacing!

%add sections	
%\let\clearpage\relax
\section{Introduction} \label{Introduction}
%Immunization effective public health intervention
Immunization is an effective tool in public health intervention for reducing child morbidity and mortality \citep{Baerninghausen2014}. In 1974, to make vaccines available to all children worldwide regardless of their social or economic status, the WHO initiated the Expanded Program on Immunization (EPI) \citep{GlobalVaccine}. Many organizations such as the WHO, the United Nations Children's Fund (UNICEF), the Global Alliance of Vaccines and Immunization (GAVI), the World Bank, national governments, and many non-governmental organizations collaborate to extend the scope and reach of the EPI. As a result of these efforts, several hundred million children worldwide have been vaccinated and millions of lives have been saved \citep{WHO2018, GAVI2020, WHO2019}. \par
%shortcomings in LMI countries
However, despite all these successes, the EPI's target of vaccinating 90\% of all children in every country by 2020 has not been achieved \citep{GlobalVaccine, WHO2019b}. Recent reports indicate that especially low- income and low-middle-income countries (LICs and LMICs) are falling short of meeting these goals. For example, out of 19.4 million children who remained un- or under-vaccinated in 2018, about 8.5 million lived in Africa \citep{WHO2019b}. Sub-Saharan Africa has a mortality rate of 77 per 1,000 live births of children under 5 years old and has a risk of child mortality 15 times higher than developed countries \citep{WHO2019b}. \par
%Reasons for these shortcomings
Several factors led to insufficient vaccine availability in LICs and LMICs. Most vaccines require a cold chain that keeps them properly refrigerated or frozen at all times until their use \citep{Vaccines2}. However, the available cold chain capacity in LICs and LMICs is limited, frequently leading to stock-outs \citep{Assi2011, cMYP2016}. Another challenge is demand uncertainty, which complicates planning, leads to insufficient vaccine availability, and high vaccine wastage. Geographically, many areas lack access to immunization services \citep{Metcalf2014, Blanford2012}. In rural areas of many African countries, vaccination coverage is negatively impacted by the limited access to healthcare services, which makes it difficult for patients to seek routine immunization services and has led to a rural-urban disparity in vaccination rates \citep{Metcalf2014, Blanford2012}. To reach children in remote areas, many countries rely on outreach activities, in which health workers travel to under-served communities to offer immunization \citep{SOS, RED}. The success of these activities depends on adequate planning and on vaccine and personnel availability. However, reaching remote areas is often challenging due inadequate road networks, long distances, and the need to maintain the cold chain's integrity during transportation. \par
%drones for the rescue
Utilizing drones for vaccine delivery may address some of the existing vaccine distribution challenges in LICs and LMICs. Recent developments in drone and battery technology allow for the economical use of drones for civil applications \citep{Zipline2019}. Drones can fly on a direct path at high speed and could enable on-demand deliveries to mitigate stock-outs caused by lacking cold chain capacities. Furthermore, the drone's independence of road infrastructure could allow it to reach remote areas. This would help communities that have been excluded from the benefits of vaccination due to lack of accessibility to medical services.\par
%existing projects
There is precedence in the use of drone for the delivery of medical supplies. In 2016, Rwanda successfully implemented drone delivery of blood for transfusions. Blood, like vaccines, also must be stored in a cooled environment, and its demand is uncertain \citep{Ackerman2019}. The implementation of drone delivery allowed centralizing the inventory, which led to a reduction of stationary storage capacity and stock-outs. In 2017, this successful project encouraged Tanzania to start using drones for delivering medical supplies. In 2018, the government of Ghana, in collaboration with UNICEF and a private partner,  built a drone delivery network for vaccines \citep{GAVI2019}. In 2018, drone-delivered vaccines reached remote communities in Vanuatu \citep{UNICEF2018}. Since 2019, through commercial contractors, Vanuatu has frequently used drones to deliver vaccines over the country.\par
%Gap in literature
Despite a growing interest in the operations research community to model vaccine supply chains (VSC) and the successful implementation  of drone delivery networks, research in this area is limited \citep{DeBoeck2020a}. Studies by \citet{Haidari2016} and \citet{Walia2018} conclude that using drones, instead of traditional land transportation, yields potential cost savings and increases vaccine availability at clinics. %\cite{Haidari2016, Walia2018}. 
However, these studies focused mainly on costs and did not consider the potential to improve vaccination rates by serving remote areas lacking access to immunization services. 
 
%Additional challenge due to covid disruptions
Drone deliveries might help save the hard-fought gains in immunization coverage in LICs and LMICs that were put at risk due to the COVID-19 pandemic \citep{WHO2020b}. By 2020, more than 13 million children missed routine immunizations, and close to 78 million children will miss routine immunizations due to COVID-19 pandemic-induced disruptions \citep{Berkley2020}. The WHO reported that in 2021 alone more than 25 million children missed out one or more doses of DTP through routine immunization services \citep{WHO2022_jul15}. 
This is due to canceled flights, temporarily suspended immunization campaigns, and people's reluctance to visit immunization session because they fear of getting infected with COVID-19. These disruptions could lead to disease outbreaks (such as measles, polio, yellow fever) which would additionally burden health systems already fighting the impacts of COVID-19 \citep{Berkley2020}.  

%Contributions
This study proposes a mixed-integer optimization model that determines a drone delivery network to maximize the number of children who are vaccinated. 
This model captures ($i$) the trade-offs between the frequency of deliveries and inventory level, ($ii$) communities' accessibility to immunization services, and ($iii$) the impact of drone deliveries for outreach activities on the efficiency of immunization programs.  
We validate the model via a case study using real-life data from Niger. We provide insights that could improve supply chain's performance.

%Organization of study
This study is organized as follows. In Section \ref{SupplyChainStructure}, we describe the structure of VSCs of LICs and LMICs and its operations. In Section \ref{RelatedLiterature}, we discuss the relevant literature. 
Next, we provide a formal problem definition and mathematical model of a VSC in Section \ref{ProblemStatementAndModel}. In Section \ref{SolutionMethod}, we discuss the approach used to solve the model. This includes a pre-processing method to reduce the problem size.  Section \ref{CaseStudy} introduces the case study.   
Finally, in Section \ref{Conclusion}, we conclude the study with a summary and potential future research directions.

\section{Background on Vaccine Supply Chains} \label{SupplyChainStructure}
\begin{figure}[!h]
	\includegraphics[width=\linewidth]{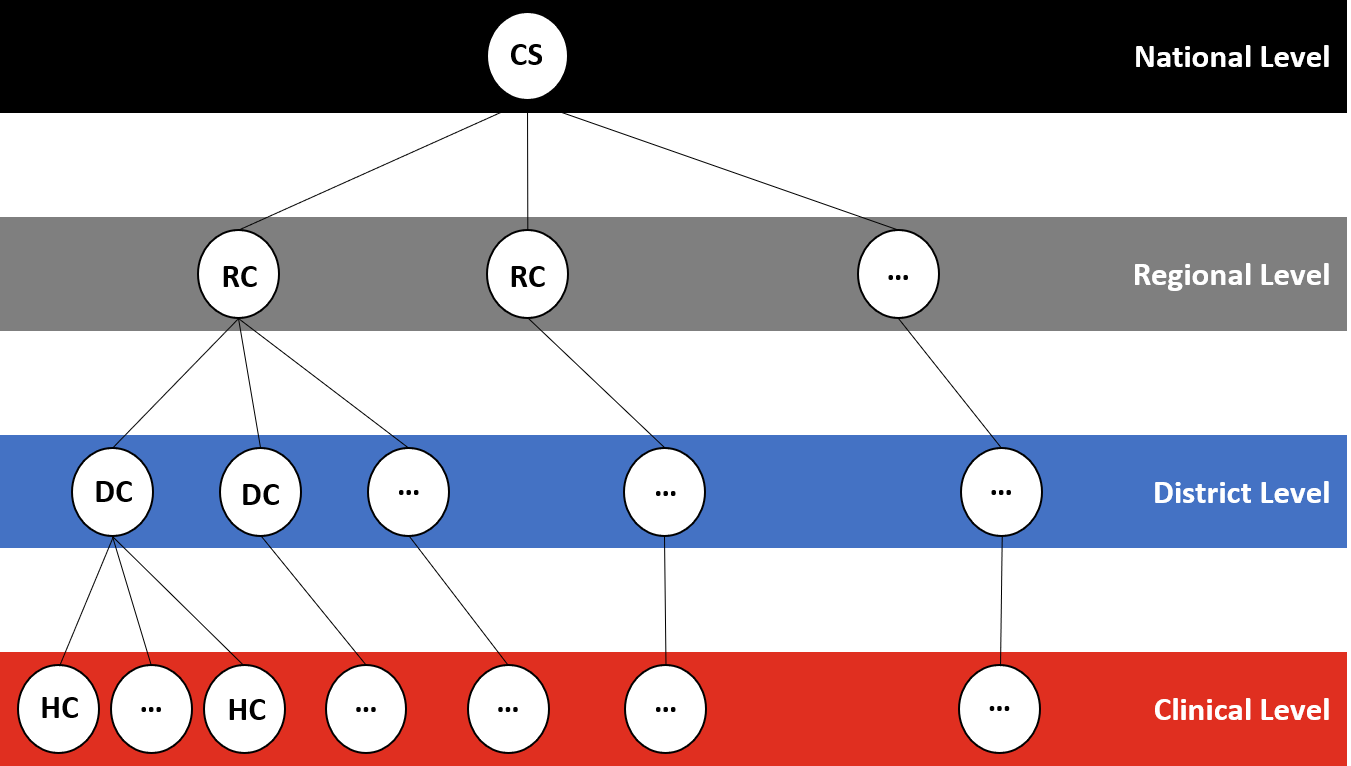}
	\caption[Typical four-level VSC]{A typical four-layer VSC (adopted from \citet{Lee2015})}
	\label{VSCStructure}
\end{figure}

Most vaccines degrade under heat exposure and must be kept at low temperatures. Few vaccines need to be kept frozen (between -25°C to - 15°C) and most need to be refrigerated (between 2°C to 8°C) \citep{Vaccines2}. Therefore, VSCs are equipped with cooling devices to maintain the required temperature ranges during storage and transportation. In LICs and LMICs, VSCs typically consist of three to five layers (i.e., echelons). In most countries,  VSCs have four layers (see Figure \ref{VSCStructure}) consisting of national-, regional-, district-level storage facilities, and clinics \citep{Lee2015}. At the national level, a central store delivers  vaccines downstream to  
regional centers (RCs) and district centers, which function as transshipment nodes. These centers deliver vaccines to clinics and health centers (HCs),  which offer immunization services to children. 

In many LICs and LMICs, major parts of the population reside in rural areas lacking access to HCs because of poor road infrastructure, long travel times, and missing or not affordable public transportation \citep{Metcalf2014, Blanford2012, Jani2014, Toik2010}. Hence, additional immunization services are often provided via outreach activities \citep{SOS, RED}. To support outreach activities, health care workers travel to rural areas.   
This service is typically provided at regular time intervals (e.g., quarterly) and pre-determined locations \citep{SOS, Lim2016}. However, the success of outreach activities depends on the availability of resources. 

Vaccines come in different vial sizes, have different cold storage requirements, and might require a diluent for reconstitution before the administration \citep{Vaccines,Vaccines2}.
World Health Organization reports that vaccines are lost, damaged, destroyed, or discarded. Vaccine wastage is typically classified into closed vial wastage (CVW) and open vial wastage (OVW). CVW primarily includes wastage due to ineffective temperature control, expiry, and breakage during handling processes. OVW occurs when multi-dose vials are opened but not completely used within a given time frame and the remaining doses have to be discarded. The OVW rate is determined by the session size, vial size, and discard time. Historically, OVW wastage contributes more to the overall wastage rate and can be quite high, with rates over 50\% \citep{WastageCalc}.

\section{Literature Review}\label{RelatedLiterature}
Relevant to our work is the literature on supply chain network design (SCND) for vaccines. \citet{Lemmens2016} surveyed the SCND literature and assessed its application to VSCs. The authors conclude that general SCND models are missing key issues unique to VSCs. Thus, researchers have to go beyond the existing SCND literature to optimize VSCs. In fact, many researches developed models explicitly for VSCs \citep{Assi2011, Azadi2020, Chen2014, de2020vaccine, Haidari2016, Lim2016, McCoy2014, Walia2018, Yang2020}; for an extensive review we refer to \citep{Duijzer2018, DeBoeck2020a}. \par

Discrete-event simulation models for VSCs were presented by \cite{Assi2011} and \cite{de2020vaccine}.  These simulation models consider cold chain capacities (inventory and transportation), shipping policies, vaccine characteristics and stochastic demand for vaccines. \cite{de2020vaccine} used their model to assess the impact of disruptions caused by the rainy season in Madagascar. The simulation model from \cite{Assi2011}, called HERMES, was used in a number of studies to evaluate the impact of different changes to the VSC such as removing a supply chain layer, altering shipping policies, or using different vial sizes \citep{Assi2011, Assi2013, Brown2014, Haidari2015, Lee2015, Lee2016, Wedlock2019}. 
However, the most closely related simulation-based study is by \cite{Haidari2016} who extended HERMES by considering drone delivery from regional and district stores to clinics, instead of using land transportation. Through a case study for Mozambique's VSC, \cite{Haidari2016} found that drone delivery can increase vaccine availability at clinics and decrease logistic costs in a broad range of settings (e.g., for different drone payloads and demand levels).\par

%add paper of responst teams and vaccine deliveries by drones
A number of optimization models have been proposed in the literature to optimize the performance of the VSCs \citep{Azadi2020, Chen2014, Lim2016, McCoy2014, Walia2018, Yang2020}. 
\cite{Yang2020} proposed a mixed-integer program (MIP) to reallocate cold capacities and reassign transportation routes to minimize inventory and transportation costs. \cite{Chen2014} proposed a linear program (LP) to make inventory, transportation, and vaccine administration decisions that maximize vaccination rates. Furthermore, the authors presented an extension to determine where in the supply chain, cold chain capacities must be added to satisfy demand at the minimum investment costs. \cite{Azadi2020} extended this model by reformulating it as a stochastic program with uncertain demand and chance constraint to ensure that demand is satisfied with a certain probability. Both, \cite{Chen2014} and \cite{Azadi2020} used their model to investigate changes to the supply chain, such as removing a supply chain layer, or changing  vaccine presentation. \cite{Walia2018} formulated an integer program (IP) to determine a cost efficient transportation mode (land transportation or drone delivery) for different transportation routes in a VSC. \cite{Lim2016} and \cite{McCoy2014} propose IPs to optimize the outcome of outreach services by improving accessibility to immunization services. \cite{McCoy2014} developed a model to determine allocation of motorcycles, needed for outreach services, to clinics to improve their services. \cite{Lim2016} compared different set covering models to maximize the number of people that have access to immunization services, assuming that only a limited number of outreach sessions can be conducted. For this purpose, \cite{Lim2016} modeled the access to services as a function of the distance between the vaccination center and locations.

\begin{table}[]
\caption{Overview of related studies}
\centering
\hfill
\begin{adjustbox}{angle=90}
\footnotesize
\begin{tabular}{|P{3.5cm}|P{3cm}|P{1.5cm}|P{2.25cm}|P{2cm}|P{1.7cm}|P{1.7cm}|P{1.7cm}|P{1.7cm}|}
\hline
{\multirow{3}{*}{\bf Publication}}                                                             & {\multirow{3}{*}{\bf Model type}} & \multicolumn{7}{P{12.8cm}|}{\bf Modeled   key issues}                                                                                                                                                                                                                                                                        \\ \cline{3-9} 
                                                                                        &                             & \multicolumn{2}{P{4cm}|}{\bf Limited Capacities} & {\bf Demand } & {\multirow{2}{*}{\bf Wastage}} & {\multirow{2}{*}{\bf Outreach}} & {\multirow{2}{*}{\bf Accessibility}} & {\bf Drone} \\ \cline{3-4}
                                                                                        &                            & \bf Inventory        & \bf Transportation       &   \bf   Uncertainty                               &                         &                           &                   & \bf  Delivery                   \\ \hline
\citet{Assi2011,   Assi2013, Lee2011,  Brown2014,    Haidari2015, Lee2011, Lee2016} & Discrete-event   simulation                      & x                & x                    & x                                                        & x                                             &                                                &                                                     &                                                      \\ \hline
\citet{Wedlock2019}                                                                             & Discrete-event   simulation                      & x                & x                    &                                                          & x                                             & x                                              &                                                     &                                                      \\ \hline
\citet{Haidari2016}                                                                            & Discrete-event   simulation                      & x                & x                    & x                                                        & x                                             &                                                &                                                     & x                                                    \\ \hline
\citet{de2020vaccine}                                                                            & Discrete-event   simulation                      & x                & x                    & x                                                        & x                                             &                                                &                                                     &                                                    \\ \hline
\citet{Lee2015}                                                                                & Discrete-event   simulation and MIP              & x                & x                    & x                                                        & x                                             &                                                &                                                     &                                                      \\ \hline
%\cite{dhamodharan2012determining}                                                             & Monte Carlo   Simulation                         &                  &                      & x                                                        & x                                             & x                                              &                                                     &                                                      \\ \hline
\citet{Chen2014}                                                                               & LP                                               & x                & x                    &                                                          & x                                             &                                                &                                                     &                                                      \\ \hline
\citet{Lim2019,   Yang2020}                                                                    & MIP                             & x                & x                    &                                                          &                                               &                                                &                                                     &                                                      \\ \hline
\citet{Walia2018}                                                                              & IP                                               &                  & x                    &                                                          &                                               &                                                & x                                                   & x                                                    \\ \hline
\citet{Azadi2020}                                                                              & SP                                               & x                &                      & x                                                        & x                                             &                                                &                                                     &                                                      \\ \hline
\citet{dhamodharan2011stochastic,   azadi2019developing}                                       & SP                                               &                  &                      & x                                                        & x                                             & x                                              &                                                     &                                                      \\ \hline
\citet{Lim2016}                                       & IP                                               &                  &                      &                                                         &                                              & x                                              &         x                                            &                                                      \\ \hline
\citet{McCoy2014}                                       & IP                                               &                  &          x            &                                                         &                                              & x                                              &                  x                                  &                                                      \\ \hline
%\cite{doerner2007multicriteria}                                       & IP                                               &                  &          x            &                                                         &                                              & x                                              &             x                                       &                                                      \\ \hline
\citet{Scott2017,  Kim2017}                                                                    & IP                                               &                  &                      &                                                          &                                               &                                                & x                                                   & x                                                    \\ \hline
\end{tabular}
\end{adjustbox}
\label{tableLitOverview}
\end{table}

\section {Problem Description and Model}\label{ProblemStatementAndModel}

The objective of this section is to develop a mathematical model to support the design of a VSC that uses drones to deliver vaccines in LICs and LMICs. 
For this purpose, we introduce problem assumptions, notation, and provide an MIP formulation.

\subsection{Assumptions}\label{ModelAssumptions}
{The focus of this research is on supporting sustained outreach activities for routine childhood immunization. Hence, we do not consider supplemental immunization programs, such as mass vaccination campaigns. Moreover, we focus on optimizing the VSC of a single country once it has vaccines available for distribution.} 

The following is a list of assumptions adopted in the proposed MIP model.
\noindent ($i$) Cold storage and cold transportation capacities are limited,
($ii$) there is sufficient capacity to store and  transport (via land) auxiliary supplies such as syringes,
($iii$) the lead time for delivering vaccines via land transportation, between any two centers in the country, is within one time period (i.e., a day),
($iv$) drone deliveries can be scheduled on-demand and have a negligible lead time,
($v$) drone deliveries are dedicated, single-trip deliveries (i.e., one beneficiary per delivery from the hub to the destination),
($vi$) drones and equipment used for cold storage and transportation are reliable and have negligible failure rates and causing no damage or vaccine spoilage,
($vii$) drones are shared among different hubs but remain stationed for a month before they can be relocated.
($viii$) there is no shortage of skilled personnel to administer vaccines at each location,
($ix$)  vaccine wastage due to breakage, incorrect storage temperatures, equipment failure, and open-vial wastage is captured via the average wastage rates,
($x$) vaccine vial size is fixed. Thus, we model the flow of doses in the supply chain. For each dose, we calculated the corresponding volume in the vial, 
($xi$) children are vaccinated at clinics or outreach posts that are located within a fixed distance from their homes. %and ($xii$) there is no backlogging of demand for vaccines.
    
\subsection{Problem Definition} \label{ProblemStatment}
The following subsections provide a formal definition of the  problem. A summary of the  notation used can be found in Appendix \ref{AppendixNotation}.\par

\subsubsection{Vaccines}
Let $\mathcal{L}$ be a set of vaccines used by a country in its routine immunization schedule. For each vaccine $l \in \mathcal{L}$, $a_l$ represents the number of doses necessary for completing the country's immunization regimen per child;  $q_{l}$ represents the effective volume per doses; and, $r_l$ represents the dilutant volume needed to reconstitute a lyophilized vaccine dose (for vaccines stored in powder presentation)\footnote{$r_l$ =0 for non-lyophilized vaccines}. Furthermore, vaccines must be kept at low temperatures until administered. There is limited cold capacity in the supply chain.  

\subsubsection{Vaccine Supply Chain}
Similar to \cite{Azadi2020} and \cite{Chen2014}, we represent the VSC for a country as a graph $\mathcal{G} = \big({\mathcal{J},~\mathcal{A}}\big)$, where $\mathcal{J}$ represent its nodes and $\mathcal{A}$ represents its arcs.
The set of nodes $\mathcal{J}$ consists of three  subsets, $\mathcal{I}, \mathcal{O}$, and $\mathcal{K}.$ Nodes in $\mathcal{I}$ represent facilities in different layers of the supply chain.  
\begin{figure}[!h]
	\includegraphics[width=\linewidth]{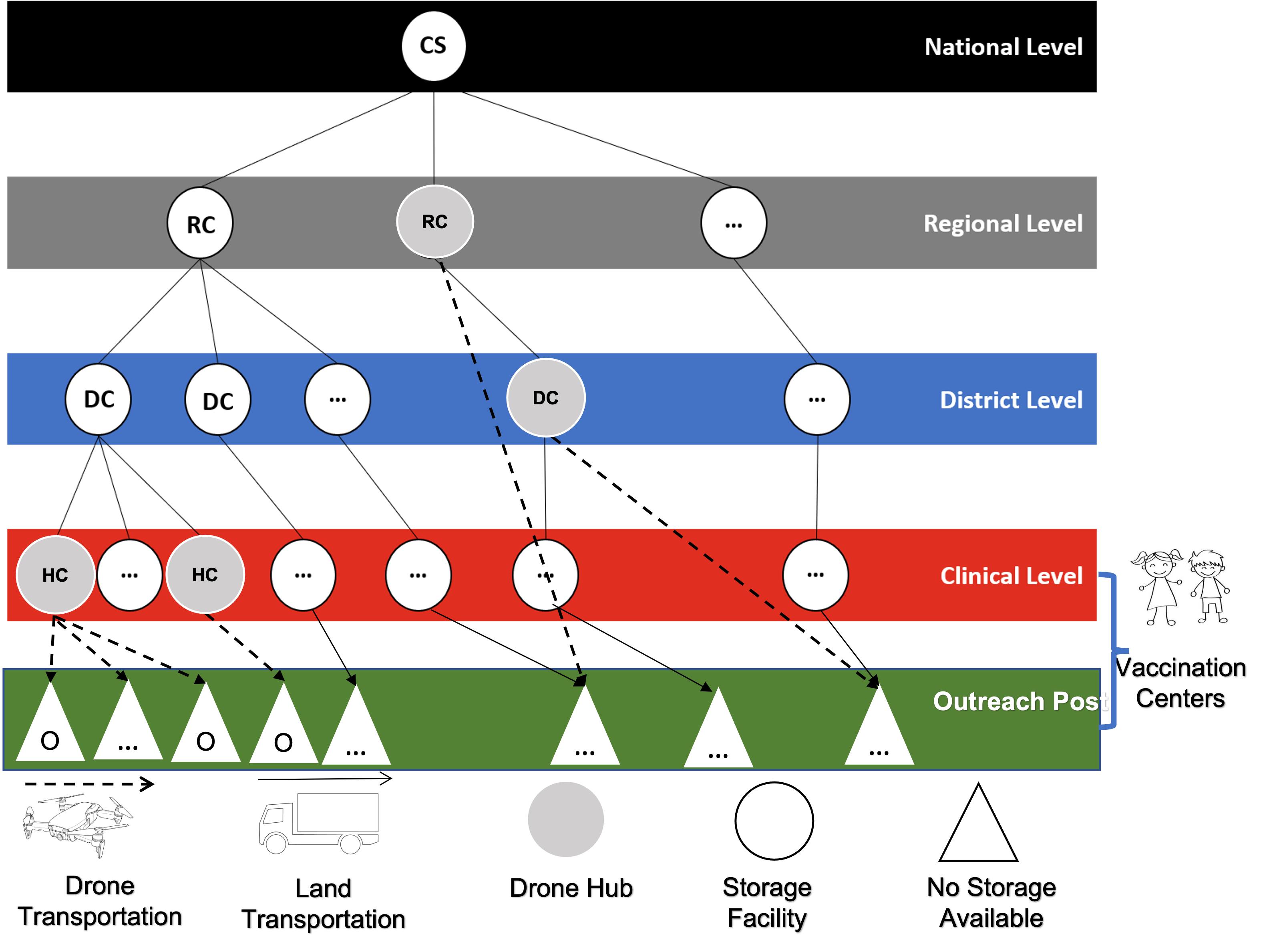}
	\caption[The proposed VSC]{The proposed VSC.}
	\label{VSCStructure2}
\end{figure}
Each  facility $i \in \mathcal{I}$ has cold storage devices with a limited holding capacity $u_{i}$. Let $w_{il}^b$ represent the expected fraction of  vaccine $l \in \mathcal{L}$ that could be lost to breakage or inadequate cooling in facility $i\in\mathcal{I}$. Nodes in $\mathcal{I}$ are potential locations for hubs since they have the available space to store and launch drones. Nodes in $\mathcal{O}$ represent the outreach posts where storage is not available. \par 

Nodes in $\mathcal{K}$ represent the communities. Let $\pi_{klt}$ represent the  demand for vaccine $l \in \mathcal{L}$ at time $t \in \mathcal{T}$ at community $k \in \mathcal{K}$. Children from these communities seek immunization services at clinics or outreach posts nearby. Let $\Pi$ represent these vaccination centers. Note that, $\mathcal{K}\cap\Pi \neq\emptyset$.

We use $w_{il}^o$ to represent open-vial wastage (OVW) during  administration of vaccine $l$ in $i \in \Pi$.   \par

The set of arcs $\mathcal{A}$ consists of three disjoint subsets; where $\mathcal{A}^\ast$ represents land transportation routes, $\dot{\mathcal{A}}$ represents potential drone transportation routes, and $\mathcal{A}^\prime$ represents walking paths from a community to a  vaccination center (clinic or outreach post). 
There is an arc between any two locations $i, j$ (thus, $(i,j) \in \dot{\mathcal{A}}$) if the distance $d_{ij}$ between them is smaller than the maximum drone range. Arcs $\mathcal{A}^\prime$ between  nodes in $i \in \Pi$ and  communities $j \in \mathcal{K}$ exist if the community $k$ is within the maximum acceptable travel distance from $i$. 
Parameters $w^t_{ijl}$ ($w^d_{ijl}$) represent the fraction of vaccines $l$  wasted during transportation between $i$ and $j$ for $(i,j) \in \mathcal{A}^\ast$ ($(i,j) \in \dot{\mathcal{A}}$) because of breakage. Parameters $u^t_{ij}$ ($u^d$) represent truck (drone) capacity along $(i,j) \in \mathcal{A}^\ast$ ($(i,j) \in \dot{\mathcal{A}}$).  
 
Figure \ref{VSCStructure2} provides a graphical representation of the proposed VSC. For simplicity, the figure does not include the communities and the walking paths to vaccination centers. The main differences between the proposed VSC and an existing/typical VSC (shown in Figure \ref{VSCStructure}) are the  addition of outreach posts, drone hubs, and drone transportation arcs.

\subsubsection{Decision Variables and Additional Notation}
Let $\mathcal{T}$ represent the planning horizon for our model. We assume that, at the beginning of the planning horizon, a centralized decision maker determines where to open drone hubs and how many drones to use. The binary variable $Y_i$ takes the value 1 if a drone hub opens in $i \in \mathcal{I}$, and takes the value 0 otherwise. The cost of opening a hub in $i$ is $c_i$. The integer variable $Z$ represents the number of drones procured. These drones are shared among hubs. The integer variable $V_{it}$ represents the number of drones used by hub $i$ in period $t$. The cost of procuring a drone with an average speed $s$ and payload $u^d$ is $\bar{c}$. The decision maker has a budget of \$$b$ for opening drone hubs and purchasing drones. This study considers that the supply chain is already equipped with vehicles for land transportation.\par

At the beginning of the planning horizon, vaccines arrive at the central store. 
The decision maker determines the number of doses of vaccine $l$ that should be delivered in period $t$ from location $i \in \mathcal{I}$ to location $j \in \mathcal{J}$ using land (via $(i,j)\in \mathcal{A}^\ast$), or using drones (via $(i,j)\in \dot{\mathcal{A}}$). These decisions,  denoted by $S_{ijlt}$ (land transportation) and $D_{ijlt}$ (drone transportation),  are affected by the availability of cold storage at the destination. 
These decisions impact the stockpile level $I_{ilt}$ at each facility $i (\in \mathcal{I})$, which in return determines the number of children that could be vaccinated in period $t$. We use $X_{ijlt}$ to represent the number of vaccines  $l$ administered at location $i (\in \Pi)$ in period $t$ to vaccinate children from community $j(\in\mathcal{K}).$  
 
\subsubsection{Objective}
%Performance measures and Objective
The objective of routine childhood immunization is to maximize the number  vaccinated  children.  A number of factors impact childhood immunizations, however,  this study focuses on increasing vaccine availability at vaccination centers.  We assume that if available, a vaccine is used to vaccinate a child in need of one. Thus, increasing availability could increase the number of children that are vaccinated.

Similar to \cite{DeBoeck2020a}, we use the supply ratio (SR) as a measure of vaccine availability, and as a performance measure for the VSC.  Let SR$_{it}$ represents the proportion of children vaccinated via the proposed VSC in period $t$ at vaccination center $i\in\Pi$. Recall that,  $\pi_{klt}$ represent the number of children in need of vaccine $l$ in community $k$ in period $t$  and $X_{iklt}$ represents the number of vaccines administered at vaccination center $i$ to children from community $k$, then,
\begin{align}
    &SR_{it} = \sum_{{l\in \mathcal{L}}} \left(\frac{\sum_{k:(k,i)\in \mathcal{A}^\prime}X_{iklt}}{\sum_{k:(k,i)\in \mathcal{A}^\prime}\pi_{klt}}\right), \hspace{0.3in} \forall i\in \Pi, t \in \mathcal{T}.
\end{align}

Let $N_k$ represent an upper bound on the number of fully immunized children of community $k$. To calculate this bound we divide the number of vaccines administered in $k$ during the planning horizon, with $\alpha_l$ -- the number of doses of vaccine $l$ necessary for completing the immunization regimen.
\begin{align} 
 &N_{k} \leq \sum_{t \in \mathcal{T}}\frac{\sum_{i:(k,i)\in \mathcal{A}^\prime}X_{iklt}}{a_{l}},~&&\forall l \in \mathcal{L},~k \in \mathcal{K}.
\end{align} 
 
The objective function, $f$, of the proposed model maximizes $N_k$ and the total number of vaccines available via the VSC. The second term in the objective maximizes SR$_{it}$. We adjust the relative importance of both terms in the objective via $\epsilon \in (0,1).$   
 \begin{align}
%objective
f = ~ & max ~\sum_{k \in \mathcal{K}} N_{k} + \epsilon  \sum_{(k,i) \in \mathcal{A}^\prime}\sum_{l\in \mathcal{L}}\sum_{t \in \mathcal{T}} X_{iklt}. \label{SSObj} 
\end{align}

\subsection{A Mixed-Integer Programming Model} \label{ModelFormulation}
%Dr. Proano this is the most recent version of the model section, the date of this comment is the 12.28.2021
We propose an MIP that extends the model proposed by \citet{Chen2014} by considering the establishment and operation of a drone network for vaccine delivery in addition to traditional land transportation. Drone deliveries improve vaccine availability at the community level. 
The proposed model is an extension of the traditional facility location-transportation model. This model integrates planning (mid-term) and operational (short-term) decisions. The planning decisions focus on the establishment of a drone network. The operational decisions focus on vaccine transportation, storage, and administration. The operation of drones in the VSC depends on the drone network. Below we introduce the proposed model formulation, which we call model ($P$).  
\begin{align}
%objective
z = max~& f \label{FSObj1}\\
%budget constraint
s.t.~&\sum_{i \in {\mathcal{I}}} c_{i}Y_{i} + \bar{c} Z \leq b, \label{budget}\\
&V_{it} \leq mY_{i},~&&\forall i \in {\mathcal{I}},~ t \in \mathcal{T}, \label{hubEstablishment}\\
& \sum_{i \in {\mathcal{I}}} V_{it} \leq Z,~&&\forall~t \in \mathcal{T}. \label{hubCapacity}
\end{align}

The objective function \eqref{FSObj1} maximizes the number of fully and partially immunized children. This number is impacted by the structure of the proposed VSC, which incorporates a drone network. Constraint \eqref{budget} indicates that the available budget limits the establishment of drone hubs and the number of drones purchased. Recall that  drones are shared in this VSC to increase  utilization. Thus, constraints \eqref{hubEstablishment} allow drones to be used by facility $i$ if a hub is already established.
Here, $m$ is the maximum number of drones that could be purchased using the available budget. Constraints \eqref{hubCapacity} force the total number of drones utilized in the VSC during time $t$ to be less than or equal to the total number of drones purchased. 
\begin{eqnarray}
I_{ilt} = (1-w_{il}^b)I_{ilt-1} - \sum_{j:(i,j)\in \mathcal{A}^\ast} S_{ijlt}
 + \sum_{j:(j,i)\in \mathcal{A}^\ast}(1-w_{jil}^t){S}_{jilt-1} - \sum_{j:(i,j)\in \dot{\mathcal{A}}} D_{ijlt} +\notag \\
+ \sum_{j:(j,i)\in \dot{\mathcal{A}}}(1 - {w}^d_{jil})D_{jilt}    - \left(\frac{\sum_{k:(k,i)\in \mathcal{A}^\prime}X_{iklt}}{1-w^o_{il}}\right), \hspace{0.2in} \forall ~i\in\mathcal{I},~l\in \mathcal{L},~t\in \mathcal{T}, \label{inventoryBalance} 
\end{eqnarray}
\vspace{-0.2in}
\begin{align} 
\sum_{l \in \mathcal{L}}q_{l}\left[I_{ilt} +  \sum_{j:(j,i)\in \mathcal{A}^\ast}(1-w^t_{jil}){S}_{jilt-1} 
 + \sum_{j:(j,i)\in \dot{\mathcal{A}}}(1 - {w}^d_{jilt})D_{jilt} \right]  \leq u_{i}, &&  ~i\in \mathcal{I},~t \in \mathcal{T}. \label{storageCapacity}
\end{align}
\begin{align}
I_{il0} = 0,   \hspace{0.5in} \forall ~i \in \mathcal{I}, ~l\in {\mathcal{L}}. \label{initialInventory}
\end{align}
Constraints \eqref{inventoryBalance} - \eqref{initialInventory} are related to the inventory. Constraint \eqref{inventoryBalance} determines inventory balance for each vaccine at each facility for each period. This constraint accounts for the loss of vaccines in the supply chain due to breakage during transportation and inventory, and OVW. The constraint also accounts for the one-period lead time for land transportation. Constraints \eqref{storageCapacity} limit  the inventory based on the available cold chain capacity. Constraint \eqref{initialInventory} initializes the inventory.  

Since the outreach posts do not carry inventory, its flow (inventory) balance constraints \eqref{inventoryBalance2} are slightly different from the inventory balance constraints for facilities in \eqref{inventoryBalance}.
\begin{eqnarray}
 \sum_{j:(j,i)\in \mathcal{A}^\ast}(1-w_{jil}^t){S}_{jilt-1} 
+ \sum_{j:(j,i)\in \dot{\mathcal{A}}}(1 - {w}^d_{jil})D_{jilt} = \left(\frac{\sum_{k:(k,i)\in \mathcal{A}^\prime}X_{iklt}}{1-w^o_{il}}\right), ~\forall ~i\in\mathcal{O},~l\in \mathcal{L},~t\in \mathcal{T}. \label{inventoryBalance2} 
\end{eqnarray}

The following set of constraints, \eqref{landTransportationCapacity} - \eqref{droneCapacity}, enforce transportation restrictions. Constraints \eqref{landTransportationCapacity} ensure  that the volume of vaccines transported by land is limited by the associated vehicle capacities. The volume of vaccines delivered via   drones is limited by the total number of  operating hours of the drone fleet since drones can fly several times between locations. Thus, the right-hand-side of constraint  \eqref{droneCapacity} represents the total available operating hours for drone deliveries at facility (hub) $i$  in period $t$.
The left-hand-side of  constraint  \eqref{droneCapacity} represents the total number of operating hours that are required to deliver vaccines from hub $i$ in period $t$. The first term in the left-hand-side determines the number of trips per period, and the second term determines the number of hours per trip. 

\begin{align}
%transportation capacity constraints
%land transportation
&\sum_{l \in {\mathcal{L}}} q_{l}S_{ijlt}\leq {u^t_{ij}},~&&\forall~(i,~j) \in \mathcal{A}^\ast,~t \in \mathcal{T}, \label{landTransportationCapacity}\\
%drone capacity
&\sum_{j:(i,j) \in \dot{\mathcal{A}}} \sum_{l \in \mathcal{L}} \left(\frac{(q_l + r_l)D_{ijlt}}{u^d}\right)\left(\frac{2d_{ij}}{s} \right) \leq hV_{it},~ &&\forall  ~i \in {\mathcal{I}},~t \in \mathcal{T}. \label{droneCapacity}
\end{align}
Constraints \eqref{demand} indicate that the  number of vaccines administered is bounded by demand for vaccines. Constraints \eqref{FIC} provide an upper bound on the number of fully immunized children (i.e., children that received the complete vaccine regimen). 
\begin{align}
%demand constraints
&\sum_{t \in \mathcal{T}}\sum_{i:(k,i)\in \mathcal{A}^\prime}{X}_{iklt} \leq \sum_{t \in \mathcal{T}}\pi_{klt},~&&\forall ~k \in \mathcal{K},~l\in \mathcal{L},\label{demand} \\
&N_{k} \leq \sum_{t \in \mathcal{T}}\frac{\sum_{i:(k,i)\in \mathcal{A}^\prime} X_{iklt}}{a_{l}},~&&~\forall k \in \mathcal{K}, l \in \mathcal{L}.\label{FIC}  
\end{align}
%variable definition constraints
Finally, we introduce the non-negativity, the binary constraints, and the integer constraints.
\begin{align}
&N_{k}\geq 0 ,~&&\forall~k \in \mathcal{K},\label{N}\\
&S_{ijlt} \geq 0,~&&\forall~(i,j) \in \mathcal{A}^\ast,~l \in \mathcal{L}, ~t \in \mathcal{T},\label{S}\\
&D_{ijlt} \geq 0,~&&\forall~(i,j) \in \dot{\mathcal{A}},~l \in \mathcal{L}, ~t \in \mathcal{T},\label{D}\\
&V_{it} \geq 0,~&&\forall~i\in{\mathcal{I}},~t \in \mathcal{T},\label{Vit}\\
&X_{iklt} \geq 0,~&&\forall ~(k,i) \in \mathcal{A}^{\prime},~l \in \mathcal{L},~t \in \mathcal{T}.\label{X}\\
&I_{ilt} \geq 0,~&&\forall~i \in \mathcal{I},~l \in \mathcal{L},~t \in \mathcal{T}, \label{I}\\
&Y_{i} \in \{0,1\},~ &&\forall i \in {\mathcal{I}}\label{binary}\\
&Z \in Z^{+} \cup \{0\}\label{integer}
\end{align}

\section{Solution Methodology} \label{SolutionMethod}
\begin{figure}[ht]
    \centering
    \includegraphics[width=\textwidth]{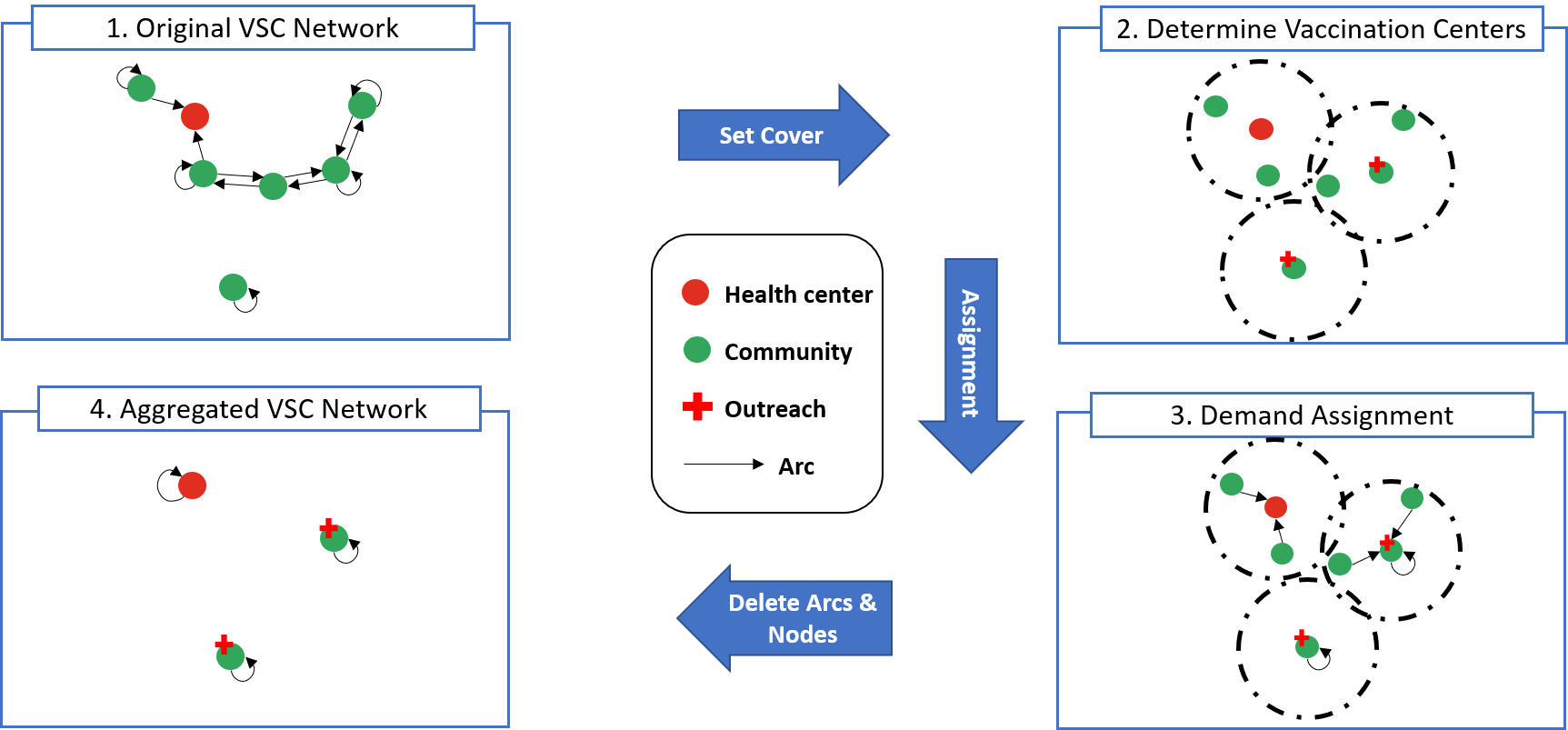}
    \caption{Pre-processing Framework}
    \label{preproc}
\end{figure}

The size of VSC, and consequently problem size, is impacted by the number of communities served. As preliminary test showed that the problem size for realistic instances is prohibitively large, we propose a pre-processing approach to reduce the  size of the VSC. In this pre-processing stage, we aggregate nearby communities into demand clusters as depicted in Figure \ref{preproc}. Thus, the problem reduces to satisfy the vaccination needs of clusters, rather than single communities. After pre-processing, the corresponding MIP model, referred as model ($Q$), is solved using the Benders Decomposition Algorithm \citep{rahmaniani2017benders}.  

The proposed pre-processing is motivated by the works of \cite{Metcalf2014, Blanford2012, Jani2014, Toik2010}  that show that people prefer getting vaccinated in nearby clinics; and that long travel times often prevent people from seeking immunization services. 
Hence, to improve  the accessibility to immunization services, we assume that communities that do not have a nearby clinic nearby must be served via outreach sessions. These outreach session can be supported by drones to ensure an adequate supply of vaccines cold chain. 
However, organizing an outreach session for every community is not efficient due to limited resources, such as, the number of available healthcare workers. Thus, it is reasonable that an outreach session pools resources to serve several communities. The challenge, however, is to determine where to locate such outreach sessions to ensure that every community has access to at least one vaccination center.
To determine what communities should serve as hosts for outreach sessions, we first solve a Set Cover model, and then use Algorithm \ref{AlgorithmDemandAssignment}, which assigns demand to the vaccination centers.

The objective of the Set Cover model is to determine the minimum number of communities that will host outreach sessions to ensure that all communities can access at least one vaccination center located within a reasonable distance (i.e., within a 5Km radius).  Consistent with the notation used in Section \ref{ProblemStatementAndModel},  $\mathcal{K}$ represents the set of communities. Let ${\bar{\mathcal{K}}}$ represent those communities that have access to a clinic ($\bar{\mathcal{K}} \subset \mathcal{K}$). Let $a_{ik}$	be an indicator parameter that takes the value 1 if vaccination center $i$ is reachable (within 5Km) from community $k$, and takes the value 0 otherwise. Let $W_{i}$ be a binary variable that takes the value 1 if a vaccination center is located in $i$, and takes the value 0 otherwise.  

\begin{align}
min~& \sum_{i \in \mathcal{K}} W_{i}, \label{SetCoverObj}\\
s.t.~& \sum_{i \in \mathcal{K}} a_{ik} W_{i} \geq 1,&& \forall~ k \in \mathcal{K}, \label{CoverEachCommunity}\\
&  W_{i} = 1,&& \forall i \in \bar{\mathcal{K}}, \label{IncludeHCs}\\
&  W_{i} \in \{0,1\}, && \forall ~ i\in \mathcal{K}. \label{SetCoverBinary}
\end{align}

The objective function \eqref{SetCoverObj}  minimizes the number of vaccination centers required to ensure access to immunization services to everyone.
Constraints \eqref{CoverEachCommunity} ensure that the outreach posts are selected such that each community has access to at least one  vaccination center within walking  distance. Furthermore, constraints \eqref{IncludeHCs} ensure that all clinics are considered as vaccination centers. \eqref{SetCoverBinary} are the binary constraints. \par

After the  vaccination centers are determined, we use Algorithm \ref{AlgorithmDemandAssignment} to assign the communities' demands to vaccination centers. Initially, each vaccination center is assigned the demand of the community where it is located. The demand of communities that are not co-located with a vaccination center but have access to at least one clinic are equally distributes among reachable clinics. For communities that have no access to a clinic, the demand is equally distributed among reachable outreach session.

\begin{algorithm}
	\caption{Demand Aggregation}
	\label{AlgorithmDemandAssignment}
	\begin{algorithmic}[1]
	        %Determine point of vaccination
		    \For{community $i \in \mathcal{K}$}
		        \If {$i$ is selected as vaccination center}
		            \State assign $i$'s demand to $i$
		        \ElsIf{$i$ has access to at least one clinic}
		            \State assign  $i$'s demand uniformly among accessible clinics
		        \Else  
		            \State assign $i$'s demand uniformly among accessible outreach posts
		        \EndIf
		    \EndFor
	\end{algorithmic}
\end{algorithm}

Algorithm \ref{AlgorithmDemandAssignment}, by aggregating demand, reduces the size of graph $\mathcal{G}$. The set of nodes in the new graph is $\mathcal{J} = \mathcal{I} \cup \mathcal{O}$, and the set of arcs is $\mathcal{A} = \mathcal{A}^\ast \cup \dot{\mathcal{A}}.$ As a result, in model ($Q$), we use decision variables $X_{ilt}$ (rather than $X_{iklt}$), which represent the total number of vaccines administered at vaccination center $i\in\Pi$. Demand for vaccinations in $i(\in\Pi)$ is $\pi_{ilt}.$ These updates impact the objective function and several constraints of ($P$). We present formulation ($Q$) in the Appendix.   
\section{Case Study} \label{CaseStudy}
In this case study, we evaluate the proposed model and solution methodology using data from Niger, 
a low-income country in West Africa, mostly covered by the Saharan desert. Roughly 80\% of its 23,196,000 inhabitants reside in rural areas with limited access to medical services \citep{cMYP2016}. Only 10\% of the roads  that  connect the main cities are paved \citep{Blanford2012}. As a result, about 60\% of the population must travel longer than 1 hour to reach a clinic. Long travel times  are one of the main reasons for low vaccination rates, which lead to immunization disparities in rural versus urban regions \citep{Blanford2012, Metcalf2014}. Therefore, road-independent drone deliveries seem promising to enhance Niger's VSC performance.

\subsection{Data Collection} \label{InputData}
Table \ref{tableDataSources} summarizes the data we have used and its sources. Additional details about the data are presented in the following subsections. 
\begin{table}[!h]
\centering
\caption{Data sources}
\label{tableDataSources}
\begin{tabular}{p{0.31\linewidth}p{0.68\linewidth}}\hline
{\bf Data}        & {\bf Source}  \\\hline
Cold chain capacities    &     \cite{cMYP2016}, \cite{Assi2011}, and \cite{Haidari2013a}\\
Population demographics           &    
\cite{census, OCHA2018, WorldBank2020, NigerPop}\\
Vaccines & \cite{Vaccines, Vaccines2, WastageCalc}               \\
Drone specifications  & \cite{Zipline2019, DoveAir2021, DHL2021, Wingcopter2021, Scott2017, Walia2018}             \\ \hline
\end{tabular}
\end{table}

%\textbf{Cold Chain:}
\subsubsection{Cold Supply Chain}
\begin{figure}[h]
    \centering
    \includegraphics[width=\textwidth]{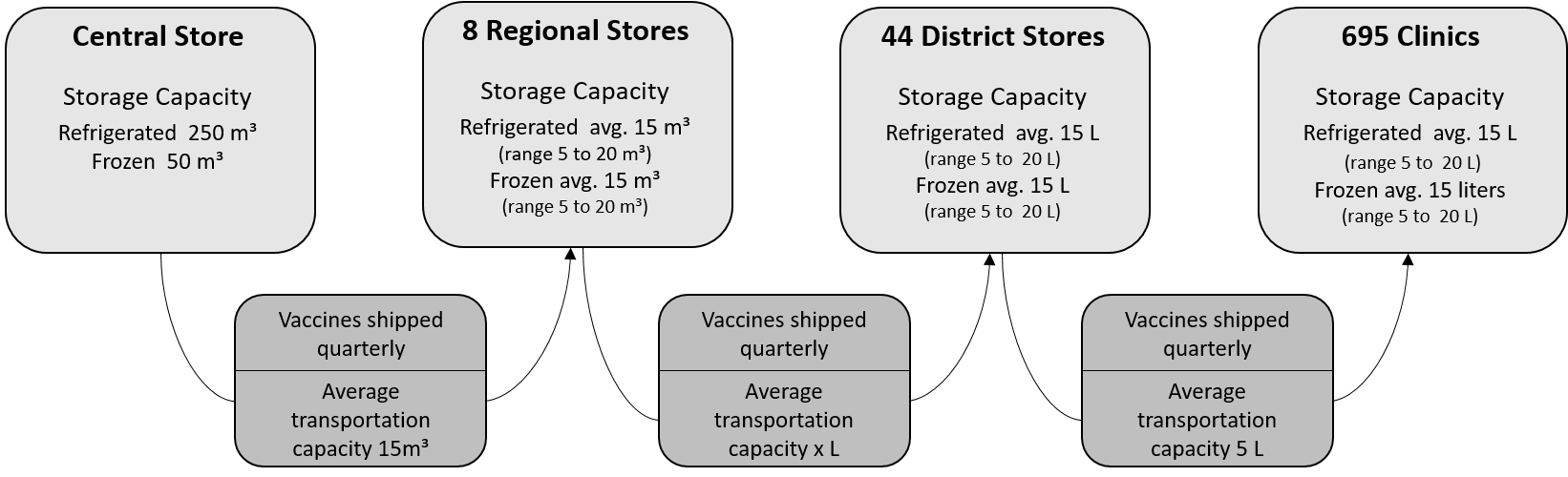}
    \caption{Niger's Cold Chain (L: liters)}
    \label{figureNigersColdChain}
\end{figure}

In Niger, vaccines are distributed through a four-tier VSC as depicted in Figure \ref{figureNigersColdChain}. The vaccines are delivered from the manufacturer to cold rooms at the central store, located in the capital, Niamey. Then, cold trucks deliver vaccines to cold rooms at eight regional stores: Agadez, Diffa, Dosso, Maradi, Niamey, Tillabery, Tahoua, and Zinder. From here, vaccines are shipped in cold boxes via small trucks to the 44 district stores, where they are stored in chest refrigerators and freezers. Finally, healthcare workers  pick up  vaccines from the district stores with cars or motorcycles, using vaccine carries for cooling, to offer immunization at 695 clinics. Currently, Niger is expanding the number of its district stores and clinics, however, we only have data for the 44 district stores and 695 clinics.
The cold capacities at clinics are  aggregated over all clinics in each region. Thus, we first allocate the regional capacity to districts based on population size, and then, allocate  districts' capacity uniformly among their clinics.
 
A clinic, on addition to servicing its patients, organizes outreach sessions to serve remote communities. Based on Niger's immunization strategy, no one should travel  longer than 5 Km to receive immunization services \citep{cMYP2016}. Since the location of the outreach posts are not published, we use the location of communities, derived from \cite{OCHA2018}, and the set cover aggregation algorithm (in Section \ref{SolutionMethod}) to determine locations of outreach posts. 
We estimated the demand for vaccination of communities using population size, and birth rates. We used the 2012 census data \citep{census} weighted by population growth rates \citep{WorldBank2020} to estimate population size of Niger in 2020 \citep{NigerPop}.
 
\subsubsection{Vaccine Characteristics}
Niger's immunization schedule for routine childhood immunization includes nine vaccines administered during six appointments as depicted in Table \ref{tableVaccines}. 
We compute vaccine wastage rates using the WHO vaccine wastage calculator \citep{WastageCalc}.

\begin{table}[!h]
	\centering
	\caption{Vaccine Characteristics}
	\small
	\label{tableVaccines}
	\begin{tabular}{P{0.175\linewidth}P{0.15\linewidth}P{0.05\linewidth}P{0.15\linewidth}P{0.15\linewidth}P{0.15\linewidth}}
		\hline
		\multirow{2}{*}{\bf Vaccine}	& \multirow{2}{*}{\bf Regimen} 				& {\bf Vial}  & {\bf Volume Per} & {\bf Volume Per}  & {\bf Storage} \\
		& & {\bf Size} &  {\bf Dose} {[}cm$^3${]} & {\bf Diluent} {[}cm$^3${]} & {\bf Requirement} \\\hline
		BCG      	& birth					              & 20        & 0.879                                 			& 0.626                                             & refrigerated or frozen \\
		DTP-HebB-Hip & 6, 10, 14 weeks			              & 10        & 3.062                                 			& -                                                 & refrigerated           \\
		Pneumo & 6, 10, 14 weeks			              & 10        & 2.109                                  			& -                                                 & refrigerated           \\
		IPV         & 14 weeks					              & 10        & 2.440                                       	& -                                                 & refrigerated           \\
		Measles     & 9, 16 months				               & 10        & 2.109                                   		& 3.142                                             & refrigerated           \\
		MenA   & 9 months					               & 10        & 2.109                                   		& 3.106                                             & refrigerated           \\
		Yellow Fever          & 9,16 months					               & 10        & 2.715                                   		& 4.730                                             & refrigerated           \\
		Rotavirus   & 6, 10 weeks				               & 1         & 17.126                                  		& -                                                 & refrigerated           \\
		OPV         & birth, 6, 10, 14 weeks	             & 1         & 0.879                                 			& -                                                 & frozen                \\\hline
	\end{tabular}
	\renewcommand{\arraystretch}{1}
\end{table}
 
\subsubsection{Drone Specifications}
%\textbf{Drone Specifications:}
A growing number of companies offer drone delivery for medical supplies, using a variety of drone technologies \citep{Zipline2019, DoveAir2021, DHL2021, Wingcopter2021, Scott2017}. Depending on the employed technology, drones have different capabilities. For instance, fixed-wing drones have larger ranges and higher speeds compared to rotor drones. However, most fixed-wing drones cannot launch vertically and need more space or special equipment for take-off and landing. Also, there is a substantial difference in terms of range and payload between battery-powered and fuel-based drones. Purchase costs for single drones and investment costs for drone hubs are hardly available since most companies offer drone delivery only as a service. Thus, we derived the drone capabilities and estimated costs using data from several drone companies \citep{Zipline2019, DoveAir2021, DHL2021, Wingcopter2021}, other academic publications \citep{ Scott2017, Walia2018}, and conversations with experts. We use  drone specifications summarized in Table \ref{tableDroneData}.

\begin{table}[!h]
\centering
\caption{Drone Specifications}
\label{tableDroneData}
\begin{tabular}{p{0.2\linewidth}p{0.245\linewidth}P{0.25\linewidth}}
\hline
\multicolumn{2}{l}{\bf Parameter} & {\bf Value} \\ \hline
\multirow{2}{*}{Range (one way)}    & battery-powered &$\sim$ 30 to 75 km\\ 
& fuel-based  & $\sim$ 900 km       \\
\multirow{2}{*}{Payload}            & battery-powered & $\sim$1.0 to 2.0 liters\\ 
&fuel-based &  \textgreater 10 liters \\
Cruise speed       &&  $\sim$ 75 km/h                                                                                                          \\
\multirow{2}{*}{Investment costs} & per drone     & \$ 30,000                                                                                                                \\
 & per drone Hub & \$ 1 Mil.                                                                                                                \\ \hline
\end{tabular}
\renewcommand{\arraystretch}{1}
\end{table}

\subsection{Experimentation and Results}
We conducted experiments to evaluate the computational performance of the pre-processing stage, and to evaluate the impact of the proposed VSC on SR. We focus on the regions Agadez, Diffa, Maradi, and Zinder and solve the model for each region separately because $(i)$ the cold capacity at the central store is not a bottleneck of the cold chain \citep{cMYP2016}, $(ii)$ each region serves its own population \citep{cMYP2016},$(iii)$ the selected regions represent diverse geographical settings, and $(iv)$ considering only a subset of regions reduces the computational complexity of the model. The experiments were run in C++ on a computer with Intel Core i7-5500 processor with 2.40 GHz using CPLEX 12.10.0 as solver. 

\subsubsection{Computational performance}
\begin{figure}[h]
     \centering
     \begin{subfigure}[b]{0.45\textwidth}
         \centering
         \includegraphics[width=\textwidth]{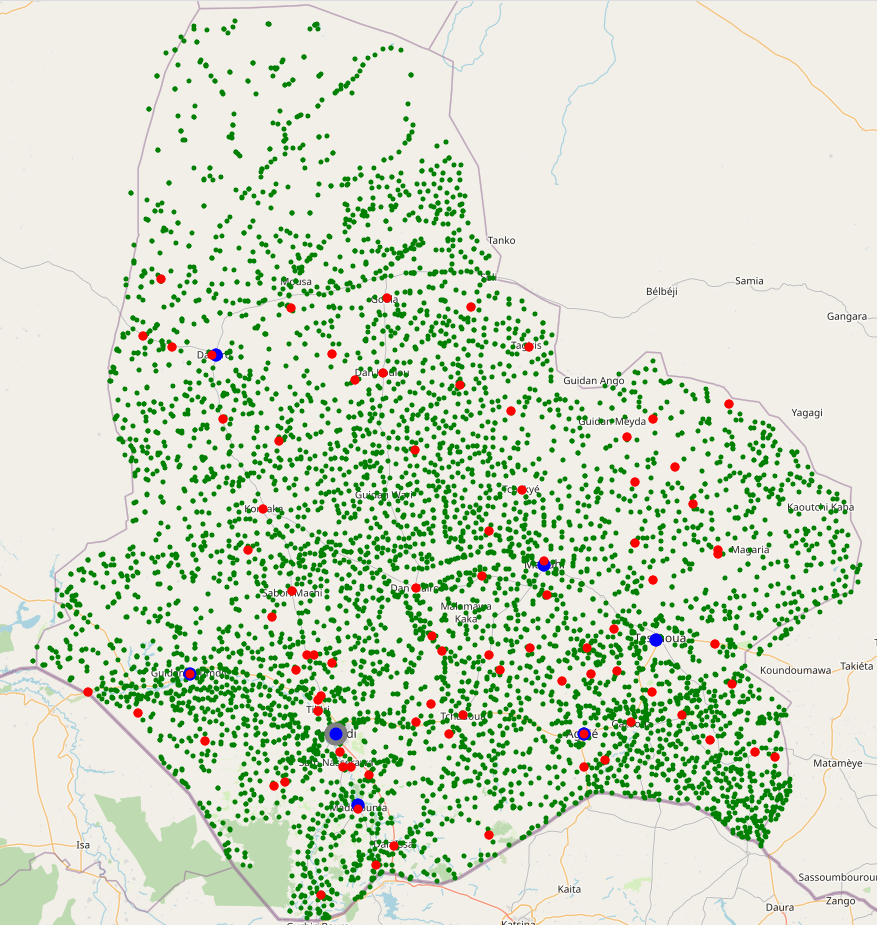}
         \caption{Original Network}
         \label{figureOriginalNetwork}
     \end{subfigure}
     \hfill
     \begin{subfigure}[b]{0.45\textwidth}
         \centering
         \includegraphics[width=\textwidth]{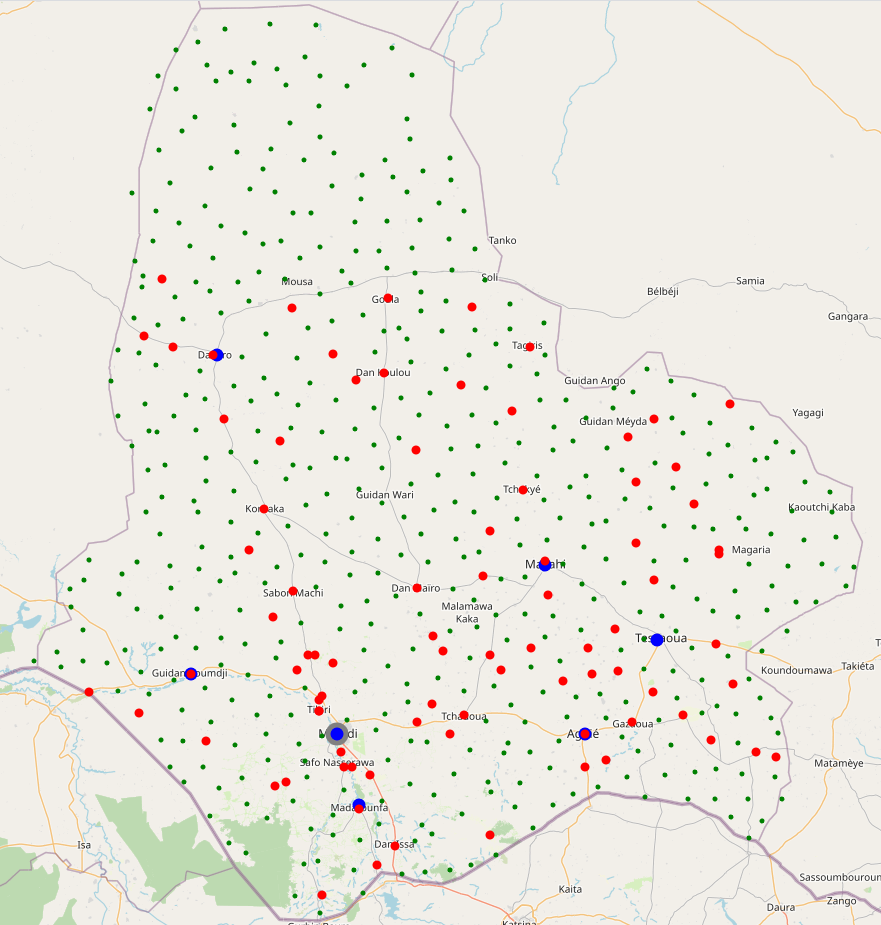}
         \caption{Aggregated Network}
         \label{figureAggregatedNetwork}
     \end{subfigure}
     \caption[Set cover aggregation applied to Maradi]{Maradi's VSC network. Gray nodes are regional centers, blue nodes are distribution centers, red nodes are clinics, and green nodes are communities.}
        \label{figureMaradiNetworkAggregation}
\end{figure}

To evaluate the impact of the pre-processing stage on the solution time and quality, we solved models ($P$) and ($Q$), which correspond to the original and aggregated networks. We used 5Km as the maximum acceptable travel distance to a vaccination center, which aligns with Niger's immunization strategy \citep{cMYP2016}. 
Table \ref{TableResultsNetworkAggregation} reports the results of the experiments. Figure \ref{figureMaradiNetworkAggregation} presents Maradi's VSC network before and after aggregation. 

\begin{table}[h]
\center
\caption[]{Network sizes and running times for the original and aggregated networks.  ($\ast$ no solution found within 10 hours)}
\label{TableResultsNetworkAggregation}
%\small
%\renewcommand{\arraystretch}{1.5}
\begin{tabular}{P{0.2\textwidth}P{0.1\textwidth}P{0.1\textwidth}P{0.1\textwidth}P{0.1\textwidth}P{0.1\textwidth}P{0.1\textwidth}} \hline
\multirow{2}{*}{\bf Region} & \multicolumn{2}{c}{\bf Nodes [-]} & \multicolumn{2}{c}{\bf Arcs [-]} & \multicolumn{2}{c}{\bf Run Time [sec]} \\
%\textbf{Region} & {\textbf{Nodes}}       & {\textbf{Arcs}}        & {\textbf{Average Run time}}  \\
                                 & \bf Orig. &\bf Agg. & \bf Orig. & \bf Agg.  & \bf Orig. & \bf Agg.  \\\hline
Agadez                           & 1283               & 634                 & 7282               & 1010                & 30    sec                       & 7    sec                                                                                                                                          \\
Diffa                            & 2144               & 629                 & 15980              & 1135                & 155 sec                           & 8   sec                                                                                                                                              \\
Maradi                           & 5383               & 569                 & 93079              & 2227                &                  $\ast$            & 81   sec                                                                                                                                                    \\
Zinder                           & 7590               & 1058                & 141262             & 2624                &               $\ast$                & 129 sec                                                                                                                       \\\hline               
\end{tabular}
\end{table}
The results show that demand aggregation significantly reduces the network sizes and running times of CPLEX. For Agadez and Diffa the average running times were reduced by 77\% and 95\%,  respectively. For Maradi's and Zinder's original network, CPLEX failed to find a solution within 10 hours. However, after the network aggregation, the instances were  solved within a few minutes. In fact, the computational challenges faced when solving the large problem instances, motivated the proposed pre-processing stage.

We found that both models, ($P$) and ($Q$), establish the same drone hubs in Agadez and Diffa. However, the SRs for model ($Q$) are 1.62\% smaller on the average.
Nevertheless, solving model ($Q$) yields near optimal solutions in reasonable time. 

\subsubsection{Discussion of results}
The following discussion addresses some relevant questions related to our proposed research. 

\noindent{\bf What is the impact of outreach sessions and drone deliveries on SR?} We establish a baseline scenario, where there are no drone deliveries and no outreach sessions. This is achieved by setting the available budget to build the drone network to zero.  Several researchers assume that families in LICs and LMICs have  access to clinics \citep{Azadi2020, Chen2014}. Additionally, these researchers assign  communities uniformly among clinics. These simplifying assumptions overestimate SR. To demonstrate this, we conducted a few experiments that solve the baseline scenario assuming full and limited access to clinics. Under the assumption of full access, we assign communities uniformly among clinics. Under the assumption of limited access, we assign a community to a clinic if located within 5Km of the clinic. Otherwise, the community does not have access to a clinic, and their children are not vaccinated. We measure the SR value for each clinic, and the SR value for each community.  

\begin{table}[!h]
\centering
\caption{Region-Based SR for the Baseline Scenario}
\label{tableBaselineRegions}
\begin{tabular}{ccccc} \hline
\multicolumn{1}{c}{\multirow{2}{*}{\textbf{Region}}} & \multicolumn{2}{c}{\textbf{Full access}}                                    & \multicolumn{2}{c}{\textbf{Limited access}}  \\
&\multicolumn{1}{c}{\textbf{Communities}} & \multicolumn{1}{c}{\textbf{Clinics}}                     &   \multicolumn{1}{c}{\textbf{Communities}} & \multicolumn{1}{c}{\textbf{Clinics}}   \\ \hline
Agadez                                      &     93.59 \%
                           &     93.59 \%                              & 41.55 \%                               &         81.83 \%                          \\
Diffa                                       &      99.51 \%                          &       99.51 \%                            &             23.05 \%                     &       88.80 \%                            \\
Maradi                                    &           70.34 \%                      &                  70.34 \%                     &       23.51 \% 
                   &          75.24 \%                          \\
Zinder                                      &            81.49 \%                    &                  81.49 \%                  &         22.19  \%                &    98.20  \%     \\ \hline      
\end{tabular}
\end{table}

Table \ref{tableBaselineRegions} summarizes the results of these experiments. Assuming full access to a clinic overestimates SR, and consequently overestimates the proportion of children that could get vaccinated. The SR values for the regions of study are lower when accessibility to clinics is considered. Based on our results, 81.83\% of the children who have access to a clinic in Agadez could be vaccinated. This is because of cold limited storage capacity at clinics. Overall, only 41.55\% of the children of Agadez could be vaccinated. This is because of most children do not have access to a clinic and/or limited cold storage capacity.

The SR values in Table \ref{tableBaselineRegions} are below 100\% when we assume full access to clinics. This is because of limited availability of cold storage and transportation in the VSC. The SR values of a region are different under the assumption of full and limited accessibility. This is because of the approach followed in assigning communities to clinics. When we assume limited access, the large population of Agadez city, Arlit, Diffa city and Maradi city is assigned to clinics located in these  cities. This leads to a disproportional number of children assigned to clinics in major cities as compared to other areas of the country, and lower SR values because of limited cold storage. Such realities could prevent Niger from achieving the EPI target of vaccinating 90\%  of the children. Thus, increasing  access to immunization services is crucial to close the gap in immunization coverage, and confirms the findings from \cite{Blanford2012}. Using outreach sessions increases accessibility and shows a great potential to improving vaccination rates. Based on these experiments, providing access to  vaccination (.i.e., via outreach sessions) could increase SR between 1.2 to 4.3 times.   

\begin{table}[]
\footnotesize
\caption{Region-Based SR Values for the Baseline Scenario: Changing Cold Capacity}
\label{tableBaselineRelievedBottlenecks}
\begin{tabular}{cccccc}\hline
\multirow{2}{*}{\textbf{Region}} & \multirow{2}{*}{\textbf{Bottleneck}} & \multicolumn{2}{c}{\textbf{Limited Capacity}} & \multicolumn{2}{c}{\textbf{Unlimited Capacity at Bottleneck}} \\
&&\multicolumn{1}{c}{\textbf{Communities}} & \multicolumn{1}{c}{\textbf{Clinics}}                     &   \multicolumn{1}{c}{\textbf{Communities}} & \multicolumn{1}{c}{\textbf{Clinics}}\\\hline
Agadez                  & Agadez City, Arlit         & 41.55\%               & 81.84\%                  & 50.77\%               & 100\%                    \\
Diffa                   & Diffa City                 & 23.05\%               & 88.80\%                  & 25.96\%               & 100\%                    \\
Maradi                  & Maradi City                & 23.51\%               & 75.24\%                  & 30.99\%                 & 100\%                    \\
Zinder                  & -                          & 22.19\%               & 98.20\%                  & 22.19\%               & 98.20\%         \\ \hline  
\end{tabular}
\end{table}

Motivated by the results presented in Table \ref{tableBaselineRegions}, we conduct one  additional experiment for the baseline scenario with limited access. We resolve the baseline scenario assuming that the major cities of Niger have unlimited cold storage capacity. We calculate the SR values for the clinics and for the overall communities in the region, and present the results in Table \ref{tableBaselineRelievedBottlenecks}.  By increasing the cold storage capacity, the SR value of clinics in Agadez, Diffa and Maradi increase. These clinics have the inventory necessary to vaccinate the children who live within 5Km. This resulted to a 12-32\% increase of SR for the corresponding regions. 

Let us summarize our {\bf observations}: ($i$)  outreach sessions supported by drone distribution networks have the potential to drastically increase accessibility to vaccines; ($ii$)  increasing cold storage of clinics located in highly  populated areas can lead to greater improvements of SR. However, such an approach does not ensure an equitable distribution of resources.   

\begin{figure}
\centering
\begin{subfigure}{.47\textwidth} 
    \includegraphics[width=.95\linewidth]{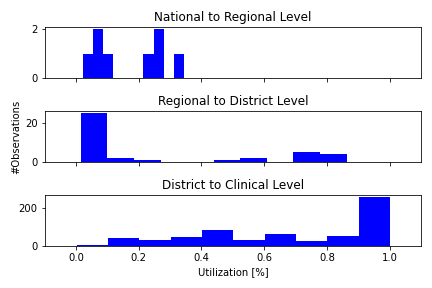}\caption{Cold Transportation} \label{fig:TCcapacity}
\end{subfigure}
\begin{subfigure}{.47\textwidth} 
    \includegraphics[width=.95\linewidth]{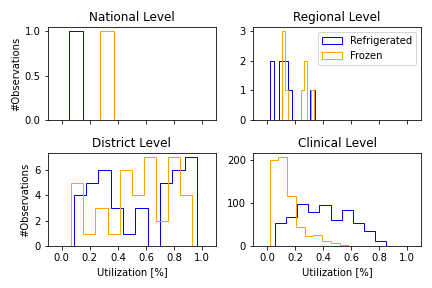}\caption{Cold Storage} \label{fig:SCcapacity}
\end{subfigure}
\caption{Resource Utilization for the Baseline Scenario } 
\end{figure}

\noindent{\bf What are the bottlenecks of the VSC?} The results presented in Table \ref{tableBaselineRelievedBottlenecks} show that limited cold storage capacity of clinics in major cities of Niger impacts SR. In Figures  \ref{fig:TCcapacity} and \ref{fig:SCcapacity} we present the utilization of cold storage and cold transportation for the baseline scenario, at all levels of the VSC. In these figures, the $x$-axis presents the utilization level, and the $y$-axis presents the number of observations. For example, trucks are used to deliver vaccines from the national-level central store to the 8 regional level stores. Since the size of vaccine vials is small, these trucks are underutilized. Cars and motorcycles are used to transport vaccines from 44 district stores to the 695 clinics. Thus, due to the small storage capacity, the utilization of these transportation vehicles is higher. Similarly, the cold storage capacity at the central store is underutilized since vaccines are moved downstream the VSC. However, for a few clinics, the cold storage utilization is over 80\%. These are clinics located in major cities.

\begin{figure}\centering
    \includegraphics[width=.9\linewidth]{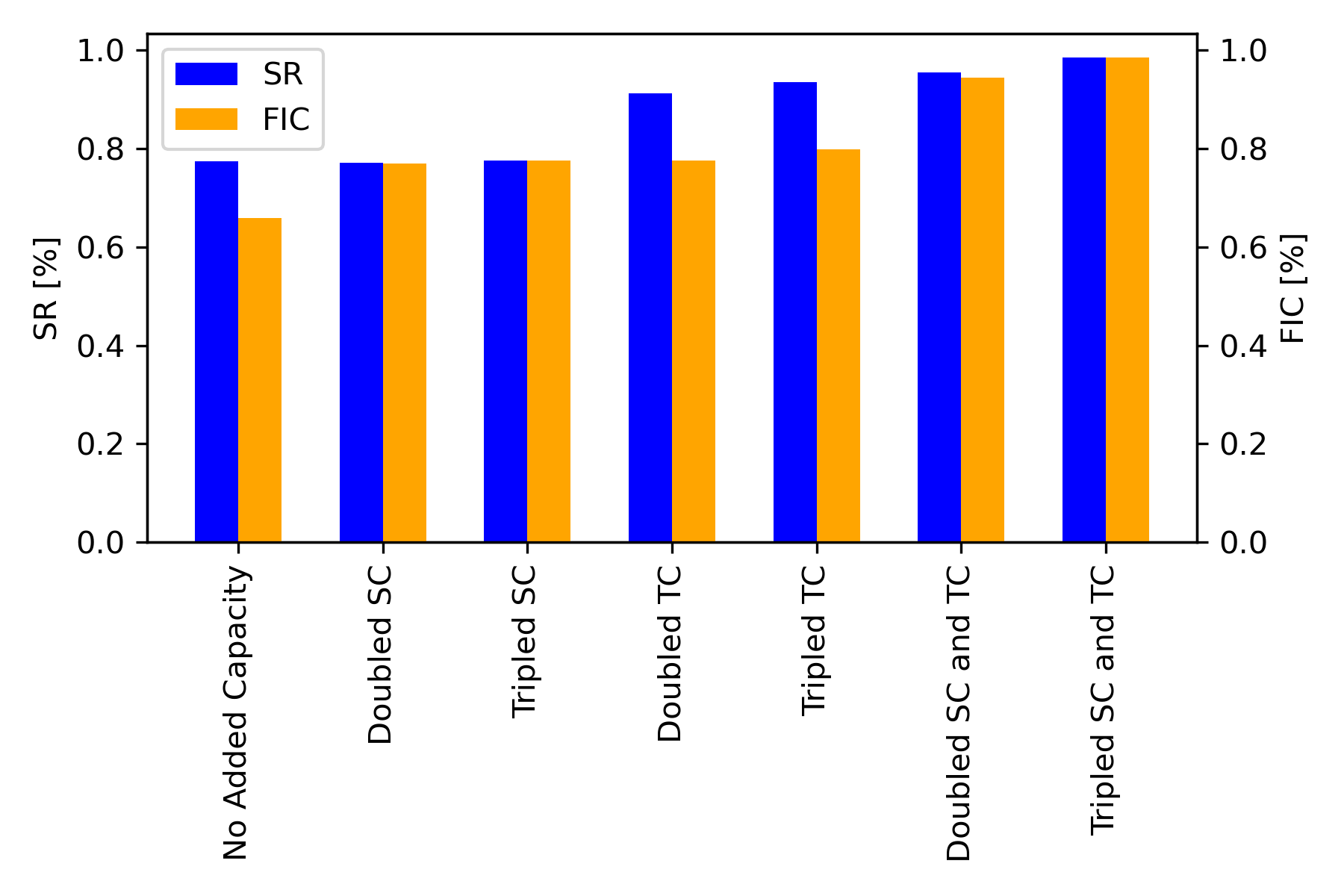}
    \label{fig:capacity}
\caption{SR \& FIC for Different Levels of Transportation and Storage Capacity of the Baseline Scenario.}\label{fig:capacity_exp}
\end{figure}

Motivated by the findings in Figures  \ref{fig:TCcapacity} and \ref{fig:SCcapacity}, we decided to evaluate the impact of increasing  storage capacity (SC) and transportation capacity (TC) on SR and the proportion of fully vaccinated children (FIC). We conducted experiments by increasing the values of TC and SC for the baseline scenario. Figure \ref{fig:capacity_exp} summarizes the results of these experiments. These are our {\bf observations}: ($i$) increasing  TC, by increasing the frequency of shipments, has a greater impact in improving SR as compared to increasing  SC; ($ii$) increasing both, TC and SC, has a greater impact in improving FIC, as compared to only increasing TC or SC.

\begin{figure}[h]
     \centering
     \begin{subfigure}[b]{0.49\textwidth}
         \centering
         \includegraphics[width=\textwidth]{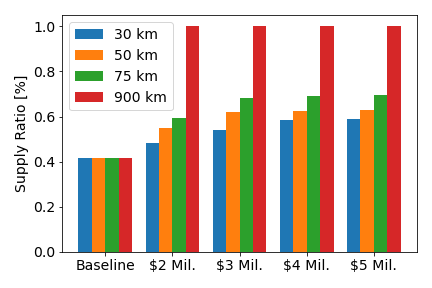}
         \caption{Agadez}
         \label{figureAgadezSensitvity}
     \end{subfigure}
     \hfill
     \begin{subfigure}[b]{0.49\textwidth}
         \centering
         \includegraphics[width=\textwidth]{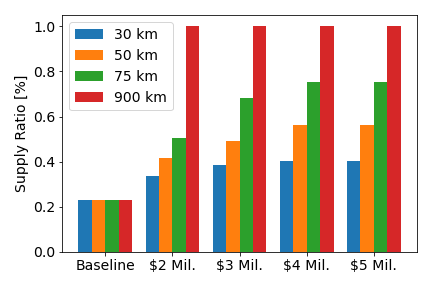}
        \caption{Diffa}
         \label{figureDiffaSensitvity}
     \end{subfigure}
     \hfill
     \centering
     \begin{subfigure}[b]{0.49\textwidth}
         \centering
         \includegraphics[width=\textwidth]{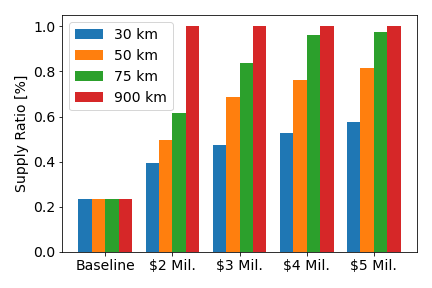}
         \caption{Maradi}
         \label{figureMaradiSensitvity}
     \end{subfigure}
      \hfill
     \begin{subfigure}[b]{0.49\textwidth}
         \centering
         \includegraphics[width=\textwidth]{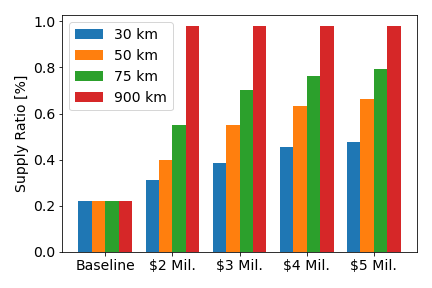}
         \caption{Zinder}
         \label{figureZinderSensitvity}
     \end{subfigure}
            \caption{Supply Ratios at Different Budget Levels and for Drone Ranges }
        \label{figureFullFactorialResults}
\end{figure}

\noindent{\bf How does budget size impact SR?} We consider  five budget levels, \$0 (baseline scenario), \$2 Millions, \$3 Millions, \$4 Millions, and \$5 Millions. These budgets are used to purchase drones with ranges 30 km, 50 km, 75 km, and 900 km. We experimented with every combination of budget and drone ranges. The  results of these experiments are summarized in Figure \ref{figureFullFactorialResults}. These are our {\bf observations}: ($i$) every positive budget results in higher SR than the baseline scenario for every region of Niger. This is because the budget is used to establish a drone delivery network; ($ii$) using drones of any range improves SR of every region. Using drones of range 900 km resulted in the highest values of SR, 100\%. However, Maradi region benefits greatly from using drones of range 50 and 75 km. This region has the smallest area, and highest population density in the country (see Table \ref{tableGeographics} in the Appendix).         

\begin{figure}[h]
	\includegraphics[width=\linewidth]{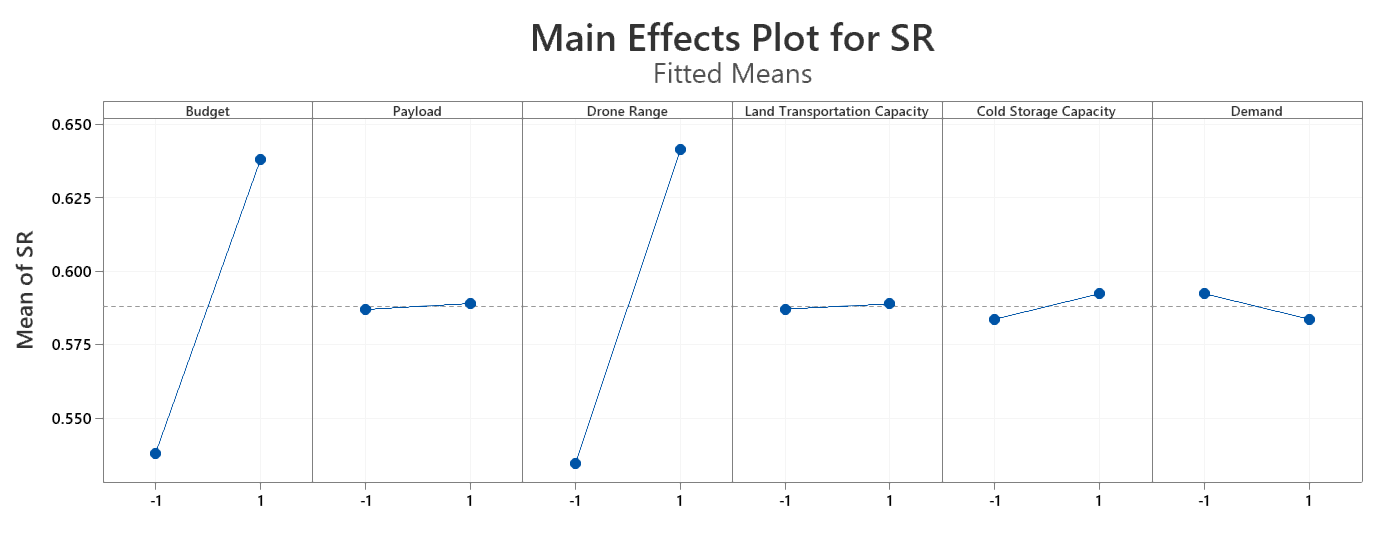}
	\caption{Plots from $2^{6-2}$ DOE Factorial Design for Agadez}
	\label{tableFactorialPlots}
\end{figure}

\noindent{\bf Do drones with higher payload or higher range have a greater impact on SR?} We initially conduct a $2^{6-2}$ factorial design of experiments (DOE) \citep{DOE} to determine what factors impact the value of SR. For each factor we determine two levels, low and high.  The factors and the corresponding levels are: budget, \$2 and \$4 million; drone payload, 1 and 2 litters; drone range, 30 and 75 km; land TC, 100\% an 125\% of the existing capacity; cold SC, 100\% an 125\% of the existing capacity; and demand, 100\% an 125\% of the existing demand. The results of these experiments for Agadez are presented in Figure \ref{tableFactorialPlots}. We make the following {\bf observations}: ($i$) drone range, rather than drone payload, greatly impacts the mean value of SR; ($ii$) budget size and drone range greatly impact the mean  value of SR; ($iii$) investments on establishing a drone network result in higher SR than investments on land TC or cold SC. 

\begin{figure}[h]
\centering
	\includegraphics[width=\linewidth]{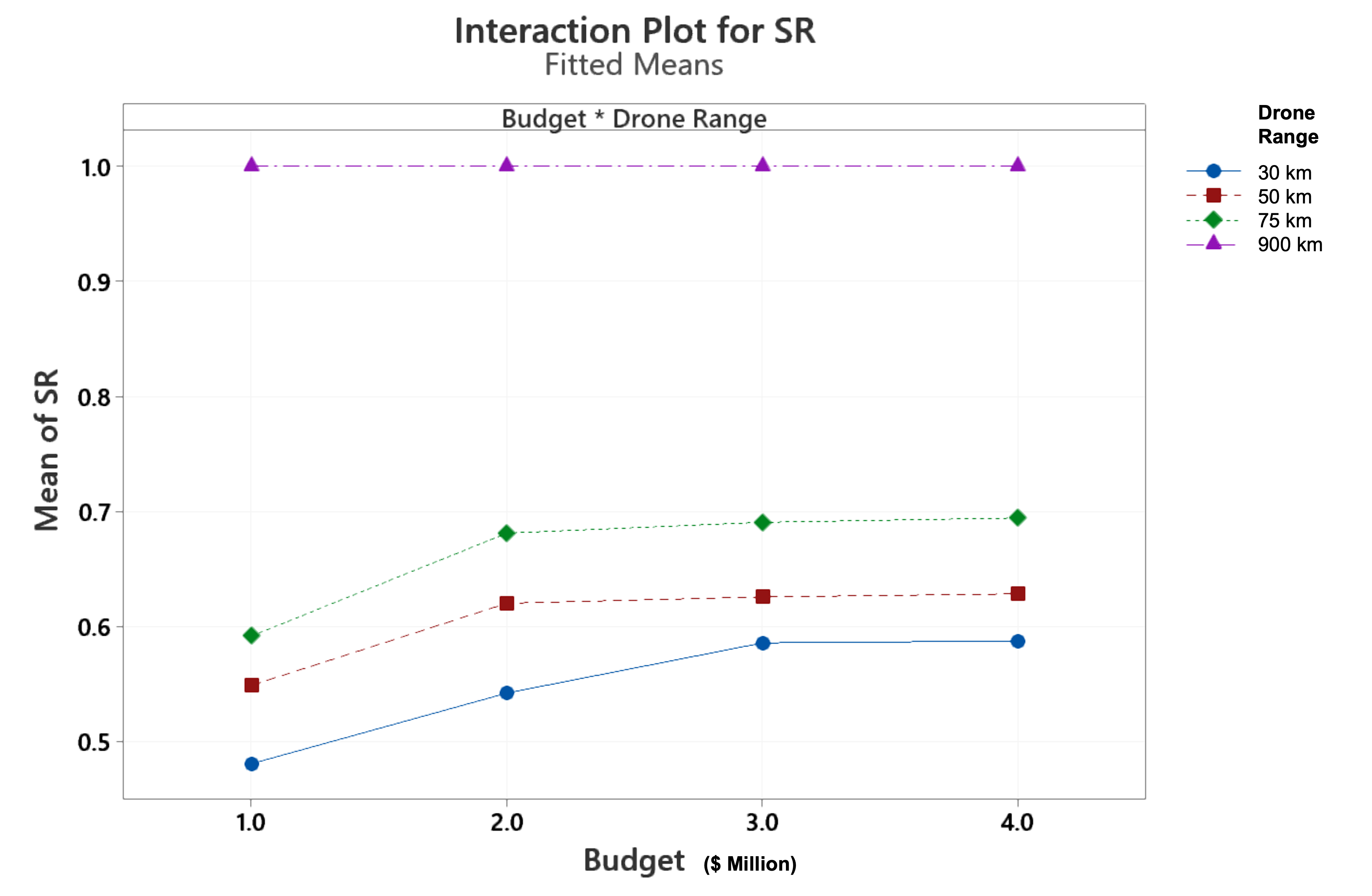}
	\caption{Results of the Full Factorial Design for Agadez}
	\label{figureInteractionAgadez}
\end{figure}

\begin{figure}[!h]
    \centering
	\includegraphics[width=0.75\linewidth]{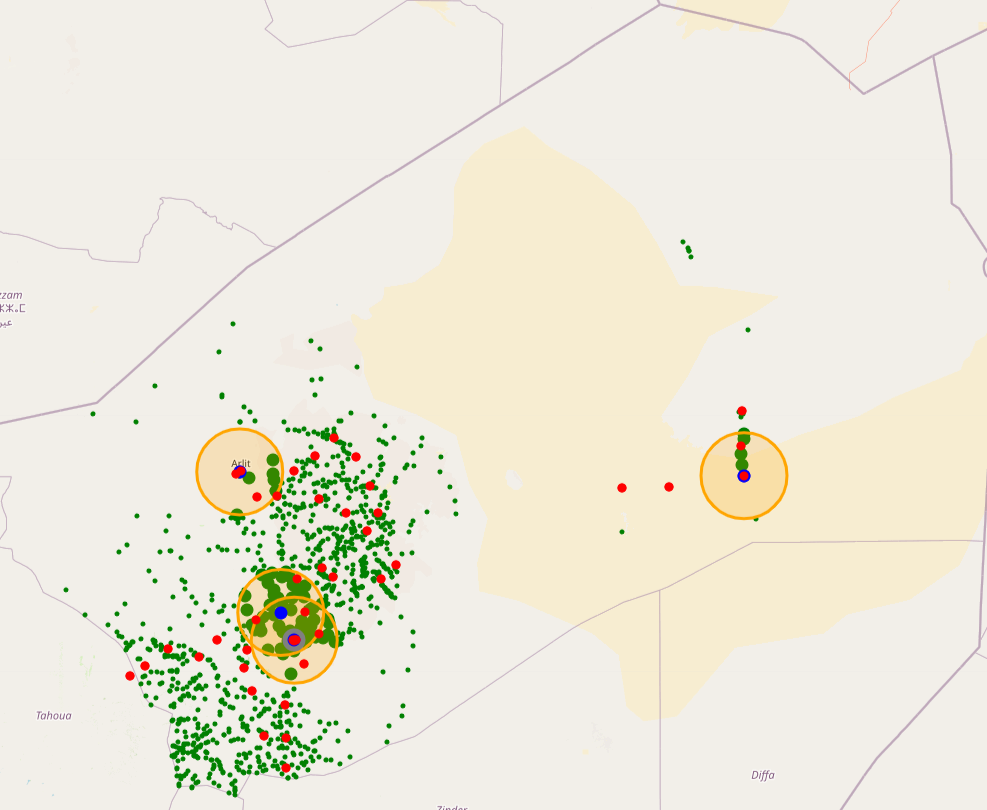}
	\caption{Drone Network of Agadez: Budget \$5 Million, Drone Range  75 km }
	\label{figureSpanAgadez}
\end{figure}

We also conducted a full factorial design to gain further understanding on the impact that budget size and drone range have on SR. Figure \ref{figureInteractionAgadez} summarizes these results for Agadez.  We make these observations: ($i$) when the budget is used to purchase short range drones, then, the initial increase of budget size greatly impacts SR. However, the maximum achievable SR from investments on short range drones is only 70\%; ($ii$) the maximum increase of SR comes from using long range drones, which improves accessibility of remote communities. 
Figure \ref{figureSpanAgadez} presents the location and coverage of drone hubs that use drones with range 75 km. The total budget to establish this drone network is \$5 million, which is used to establish 4 hubs and purchase the drones. The range of the drones does not cover communities located in remote areas.   
 
\noindent{\bf What strategies to use in designing the drone network?} Based on the results of Figure \ref{figureFullFactorialResults}, budget size impacts decisions about the number and location of drone hubs, and the fleet size. In these experiments, we observed that drone hubs are typically located at regional centers. The inventory level of regional centers is larger than district centers and clinics. Drone deliveries  provide clinics and outreach posts with access to the inventory of regional centers. 

Since LICs and LMICs have limited resources, it is of interest to evaluate what strategies these countries can use to establish and expand the drone network to increase accessibility to vaccines. Our \emph{proposed strategy} is to use regional centers as potential location for drone hubs. If budget is limited, priority should be given to regional centers that are close to highly populated areas. As additional budget becomes available, this network can be extended.  

In order to evaluate the performance of this proposed strategy, we design the following experiment that mimics a potential expansion of the drone network. The experiment solves a variation of model ($Q$) in which only regional centers are candidate locations to establish a drone hub. Let us call this model ($\bar{Q}$). We solve model ($\bar{Q}$) the first time assuming a \$2 million budget  to establish the drone hub and purchase drones. We solve model ($\bar{Q}$) a second time assuming a budget of \$3 millions, of which \$2 million is already assigned to establishing the hub(s) open the first time we solved model ($\bar{Q}$).   We continue in this fashion by changing the budget in increments of \$1 million and resolving ($\bar{Q}$) till the total budget reaches \$5 millions. The SR values for each of the solutions found are compared to the optimal solution of ($Q$) (at the  corresponding budget level). The corresponding percentage decrease of SR values are summarized in Table \ref{tableSeqOpening}. We {\bf observe} that the proposed strategy is  optimal in 78\% of the problems solved, and in the worst case reduces SR by at most 2.12\%. 
\begin{table}[!h]
\centering
\small
\renewcommand{\arraystretch}{1}
\caption{Percentage Difference of SR of the Proposed Drone Network Expansion and the Optimal Solution}
\label{tableSeqOpening}
\begin{tabular}{|c|cccc|c|cccc|} \hline
\multicolumn{5}{|c|}{\textbf{Agadez}}                                                & \multicolumn{5}{c|}{\textbf{Diffa}}                                                   \\\hline
{Budget}   & \multicolumn{4}{c|}{Drone Range {[}km{]}} & {Budget }   & \multicolumn{4}{c|}{Drone Range   {[}km{]}} \\
   {[}\$ Mil.{]}                                     & 30       & 50       & 75       & 900     &    {[}\$ Mil.{]}                                     & 30        & 50       & 75       & 900      \\\hline
2.0                  & 0.40\%               & 1.36\%               & 2.12\%               & 0.00\%               & 2.0                  & 0.00\%               & 0.00\%               & 0.71\%               & 0.00\%               \\
3.0                  & 0.00\%               & 1.36\%               & 2.12\%               & 0.00\%               & 3.0                  & 0.00\%               & 0.00\%               & 2.10\%               & 0.00\%               \\
4.0                  & 0.00\%               & 0.00\%               & 0.00\%               & 0.00\%               & 4.0                  & 0.00\%               & 0.00\%               & 0.00\%               & 0.00\%               \\
5.0                  & 0.00\%               & 0.00\%               & 0.00\%               & 0.00\%               & 5.0                  & *                    & *     &*  &* \\\hline
\multicolumn{5}{|c|}{\textbf{Maradi}}                                                & \multicolumn{5}{c|}{\textbf{Zinder}}                                                  \\\hline
{Budget} & \multicolumn{4}{c|}{Drone Range {[}km{]}} & {Budget } & \multicolumn{4}{c|}{Drone Range   {[}km{]}} \\
    {[}\$   Mil.{]}                                     & 30       & 50       & 75       & 900     &   {[}\$   Mil.{]}                                      & 30        & 50       & 75       & 900      \\\hline
2.0                  & 0.00\%               & 0.00\%               & 0.00\%               & 0.00\%               & 2.0                  & 2.31\%               & 1.91\%               & 1.10\%               & 0.00\%               \\
3.0                  & 0.00\%               & 0.00\%               & 0.00\%               & 0.00\%               & 3.0                  & 0.78\%               & 1.16\%               & 0.00\%               & 0.00\%               \\
4.0                  & 0.16\%               & 0.00\%               & 0.00\%               & 0.00\%               & 4.0                  & 0.00\%               & 1.99\%               & 0.00\%               & 0.00\%               \\
5.0                  & 0.15\%               & 0.00\%               & 0.00\%               & 0.00\%               & 5.0                  & 0.00\%               & 1.62\%               & 0.00\%               & 0.00\%       \\\hline
\end{tabular}
\end{table}
\section{Summary of Results and Future Research Directions} \label{Conclusion}
Achieving vaccination targets in low-income and low-middle-income countries (LICs and LMICs) is a challenging task. Therefore, this study examines the use of drones to improve vaccine distribution. 
For this purpose, we develop a mixed-integer linear program  that determines where to establish drone hubs to improve access to vaccination. The model captures transportation and inventory replenishment decisions in the vaccine supply chain (VSC). The model considers the spatial distribution of communities and the location of clinics in determining locations for drone hubs. We develop a case study using data from  Niger. We develop a pre-processing algorithm to reduce the problem size. The aggregated model is solved via CPLEX. 

Via the case study, we make a number of observations. The most relevant findings are:
\begin{itemize}
    \item Improving  accessibility to immunization services via outreach sessions is crucial to achieve vaccination targets $>$90\%.
    \item Drones are capable to increase vaccination rates by mitigating bottlenecks in the cold chain via on-demand deliveries to clinics and by reaching undeserved hard-to-reach areas.
    \item The benefits of drone deliveries are mainly limited by the size of the drone network and reach of the drones. Thus, it is important to determine the location of drone hubs and the drone range  to achieve vaccination targets.
    \item Our proposed strategy of using regional centers in highly populated areas as potential location for drone hubs, and extending the drone network as the budget becomes available, leads to high SR values (high quality solutions).
\end{itemize}

The efficacy of the proposed model is limited by assumptions made to deal with missing or incomplete data. 
Our proposed model does not consider vaccine hesitancy, and assumes that if vaccines are available within reach (within 5Km from the community), they are used. However, the actual demand for vaccines is influenced by many factors such as religious beliefs, education, attitude towards vaccination, accessibility of services, birth rates, and others \cite{Toik2010}.  We did not have the data and knowledge about the impact of vaccine hesitancy on vaccination demand at the community level. The model can be extended by estimating vaccine hesitancy, and using these estimated to update demand in the VSC. Changes in demand impact the design of the drone network.

The proposed model considers dedicated, single-trip deliveries. Such a  policy is simple and  easy to implement. In contrast, supplying multiple locations in one trip could reduced the number of drones used in the VSC, and shift attention to the drone payload. We only consider the initial investment cost for the drone network but do not include the operational costs  due to lack of data. Extending the model to consider these options could improve the operations of VSC. 

We assume that vaccination demand is deterministic. However, the model can be extended to evaluate the impact of stochastic demand on VSC decisions. Additionally, the model can be extended to capture the impacts of using drones for the delivery of vaccines during VSC disruptions caused by disasters, war and conflicts, pandemics, etc. 

Finally, our case study focuses on Niger. The VSC and demographics of Niger are different from other Sub-Saharan countries. However, these countries have large rural population that remains under-served because of lacking access to clinics \cite{Metcalf2014}. Therefore, it can be assumed that in these countries, drone delivery has a similar impact and helps to improve vaccination rates by serving people in rural areas. \par

\bibliographystyle{chicago}
\spacingset{1}
\bibliography{source}

\begin{thebibliography}{}

\bibitem[\protect\citeauthoryear{??}{Glo}{2013}]{GlobalVaccine}
 (2013).
\newblock {Global Vaccine Action Plan}.
\newblock {\em Vaccine\/}~{\em 31}.

\bibitem[\protect\citeauthoryear{Ackerman and Koziol}{Ackerman and
  Koziol}{2019}]{Ackerman2019}
Ackerman, E. and M.~Koziol (2019).
\newblock The blood is here: Zipline's medical delivery drones are changing the
  game in rwanda.
\newblock {\em IEEE Spectrum\/}~{\em 56\/}(5), 24--31.

\bibitem[\protect\citeauthoryear{Alliance}{Alliance}{2019}]{GAVI2019}
Alliance, G.~V. (2019).
\newblock {Ghana launches the world's largest vaccine drone delivery network}.
\newblock Available at
  https://www.gavi.org/ghana-launches-the-world-s-largest-vaccine-drone-delivery-network,
  last accessed 02.19.2020.

\bibitem[\protect\citeauthoryear{Alliance}{Alliance}{2020}]{GAVI2020}
Alliance, T. G.~V. (2020).
\newblock {Gavi @ 20 The Story of an Alliance that Today Protects Half the
  World's Children}.
\newblock Available at https://www.gavi.org/gavi-at-20, last accessed
  02.17.2020.

\bibitem[\protect\citeauthoryear{Assi, Brown, Djibo, Norman, Rajgopal, Welling,
  Chen, Bailey, Kone, Kenea, Connor, Wateska, Jana, Wisniewski, {Van Panhuis},
  Burke, and Lee}{Assi et~al.}{2011}]{Assi2011}
Assi, T.~M., S.~T. Brown, A.~Djibo, B.~A. Norman, J.~Rajgopal, J.~S. Welling,
  S.~I. Chen, R.~R. Bailey, S.~Kone, H.~Kenea, D.~L. Connor, A.~R. Wateska,
  A.~Jana, S.~R. Wisniewski, W.~G. {Van Panhuis}, D.~S. Burke, and B.~Y. Lee
  (2011, dec).
\newblock {Impact of changing the measles vaccine vial size on Niger's vaccine
  supply chain: A computational model}.
\newblock {\em BMC Public Health\/}~{\em 11\/}(1), 425.

\bibitem[\protect\citeauthoryear{Assi, Brown, Kone, Norman, Djibo, Connor,
  Wateska, Rajgopal, Slayton, and Lee}{Assi et~al.}{2013}]{Assi2013}
Assi, T.~M., S.~T. Brown, S.~Kone, B.~A. Norman, A.~Djibo, D.~L. Connor, A.~R.
  Wateska, J.~Rajgopal, R.~B. Slayton, and B.~Y. Lee (2013, jun).
\newblock {Removing the regional level from the Niger vaccine supply chain}.
\newblock {\em Vaccine\/}~{\em 31\/}(26), 2828--2834.

\bibitem[\protect\citeauthoryear{Azadi, Eksioglu, and Geismar}{Azadi
  et~al.}{2020}]{Azadi2020}
Azadi, Z., S.~D. Eksioglu, and H.~N. Geismar (2020).
\newblock Optimization of distribution network configuration for pediatric
  vaccines using chance constraint programming.
\newblock {\em arXiv preprint arXiv:2006.05488\/}.

\bibitem[\protect\citeauthoryear{Azadi, Gangammanavar, and Eksioglu}{Azadi
  et~al.}{2019}]{azadi2019developing}
Azadi, Z., H.~Gangammanavar, and S.~Eksioglu (2019).
\newblock Developing childhood vaccine administration and inventory
  replenishment policies that minimize open vial wastage.
\newblock {\em Annals of Operations Research\/}, 1--33.

\bibitem[\protect\citeauthoryear{B{\"a}rnighausen, Berkley, Bhutta, Bishai,
  Black, Bloom, Constenla, Driessen, Edmunds, Evans, et~al.}{B{\"a}rnighausen
  et~al.}{2014}]{Baerninghausen2014}
B{\"a}rnighausen, T., S.~Berkley, Z.~A. Bhutta, D.~M. Bishai, M.~M. Black,
  D.~E. Bloom, D.~Constenla, J.~Driessen, J.~Edmunds, D.~Evans, et~al. (2014).
\newblock Reassessing the value of vaccines.
\newblock {\em The Lancet Global Health\/}~{\em 2\/}(5), e251--e252.

\bibitem[\protect\citeauthoryear{Blanford, Kumar, Luo, and MacEachren}{Blanford
  et~al.}{2012}]{Blanford2012}
Blanford, J.~I., S.~Kumar, W.~Luo, and A.~M. MacEachren (2012).
\newblock {It's a long, long walk: accessibility to hospitals, maternity and
  integrated health centers in Niger}.
\newblock {\em International Journal of Health Geographics\/}~{\em 11}, 1--15.

\bibitem[\protect\citeauthoryear{Brown, Schreiber, Cakouros, Wateska, Dicko,
  Connor, Jaillard, Mvundura, Norman, Levin, Rajgopal, Avella, Lebrun,
  Claypool, Paul, and Lee}{Brown et~al.}{2014}]{Brown2014}
Brown, S.~T., B.~Schreiber, B.~E. Cakouros, A.~R. Wateska, H.~M. Dicko, D.~L.
  Connor, P.~Jaillard, M.~Mvundura, B.~A. Norman, C.~Levin, J.~Rajgopal,
  M.~Avella, C.~Lebrun, E.~Claypool, P.~Paul, and B.~Y. Lee (2014, jul).
\newblock {The benefits of redesigning Benin's vaccine supply chain}.
\newblock {\em Vaccine\/}~{\em 32\/}(32), 4097--4103.

\bibitem[\protect\citeauthoryear{Chen, Norman, Rajgopal, Assi, Lee, and
  Brown}{Chen et~al.}{2014}]{Chen2014}
Chen, S.~I., B.~A. Norman, J.~Rajgopal, T.~M. Assi, B.~Y. Lee, and S.~T. Brown
  (2014, aug).
\newblock {A planning model for the WHO-EPI vaccine distribution network in
  developing countries}.
\newblock {\em IIE Transactions\/}~{\em 46\/}(8), 853--865.

\bibitem[\protect\citeauthoryear{De~Boeck, Decouttere, J{\'o}nasson, and
  Vandaele}{De~Boeck et~al.}{2020}]{de2020vaccine}
De~Boeck, K., C.~Decouttere, J.~O. J{\'o}nasson, and N.~Vandaele (2020).
\newblock Vaccine supply chains in resource-limited settings: Mitigating rainy
  season disruptions.
\newblock {\em Available at SSRN 3584213\/}.

\bibitem[\protect\citeauthoryear{{De Boeck}, Decouttere, and Vandaele}{{De
  Boeck} et~al.}{2020}]{DeBoeck2020a}
{De Boeck}, K., C.~Decouttere, and N.~Vandaele (2020).
\newblock {Vaccine distribution chains in low- and middle-income countries: A
  literature review}.
\newblock {\em Omega (United Kingdom)\/}~{\em 97}.

\bibitem[\protect\citeauthoryear{{Deutsche Post Group}}{{Deutsche Post
  Group}}{2021}]{DHL2021}
{Deutsche Post Group} (2021).
\newblock {DHL Parcelcopter}.
\newblock Available at
  https://www.dpdhl.com/en/media-relations/specials/dhl-parcelcopter.html, last
  accessed 08.18.2020.

\bibitem[\protect\citeauthoryear{Dhamodharan, Proano, and Kumar}{Dhamodharan
  et~al.}{2011}]{dhamodharan2011stochastic}
Dhamodharan, A., R.~A. Proano, and S.~Kumar (2011).
\newblock A stochastic approach to determine the optimal vaccine vial size.
\newblock In {\em IIE Annual Conference. Proceedings}.

\bibitem[\protect\citeauthoryear{{Dove Air}}{{Dove Air}}{2021}]{DoveAir2021}
{Dove Air} (2021).
\newblock {Technical Specifications}.
\newblock Available at https://www.doveair.org/our-doves/, last accessed
  08.18.2020.

\bibitem[\protect\citeauthoryear{Duijzer, Jaarsveld, and Dekker}{Duijzer
  et~al.}{2018}]{Duijzer2018}
Duijzer, L.~E., W.~Jaarsveld, and R.~Dekker (2018, 01).
\newblock Literature review - the vaccine supply chain.
\newblock {\em European Journal of Operational Research\/}.

\bibitem[\protect\citeauthoryear{Haidari, Brown, Ferguson, Bancroft, Spiker,
  Wilcox, Ambikapathi, Sampath, Connor, and Lee}{Haidari
  et~al.}{2016}]{Haidari2016}
Haidari, L.~A., S.~T. Brown, M.~Ferguson, E.~Bancroft, M.~Spiker, A.~Wilcox,
  R.~Ambikapathi, V.~Sampath, D.~L. Connor, and B.~Y. Lee (2016, jul).
\newblock {The economic and operational value of using drones to transport
  vaccines}.
\newblock {\em Vaccine\/}~{\em 34\/}(34).

\bibitem[\protect\citeauthoryear{Haidari, Connor, Wateska, Brown, Mueller,
  Norman, Schmitz, Paul, Rajgopal, Welling, Leonard, Chen, and Lee}{Haidari
  et~al.}{2013}]{Haidari2013a}
Haidari, L.~A., D.~L. Connor, A.~R. Wateska, S.~T. Brown, L.~E. Mueller, B.~A.
  Norman, M.~M. Schmitz, P.~Paul, J.~Rajgopal, J.~S. Welling, J.~Leonard, S.-I.
  Chen, and B.~Y. Lee (2013, may).
\newblock {Augmenting Transport versus Increasing Cold Storage to Improve
  Vaccine Supply Chains}.
\newblock {\em PLoS ONE\/}~{\em 8\/}(5).

\bibitem[\protect\citeauthoryear{Haidari, Wahl, Brown, Privor-Dumm,
  Wallman-Stokes, Gorham, Connor, Wateska, Schreiber, Dicko, Jaillard, Avella,
  and Lee}{Haidari et~al.}{2015}]{Haidari2015}
Haidari, L.~A., B.~Wahl, S.~T. Brown, L.~Privor-Dumm, C.~Wallman-Stokes,
  K.~Gorham, D.~L. Connor, A.~R. Wateska, B.~Schreiber, H.~Dicko, P.~Jaillard,
  M.~Avella, and B.~Y. Lee (2015, jun).
\newblock {One size does not fit all: The impact of primary vaccine container
  size on vaccine distribution and delivery}.
\newblock {\em Vaccine\/}~{\em 33\/}(28), 3242--3247.

\bibitem[\protect\citeauthoryear{Jani, De~Schacht, Jani, and Bjune}{Jani
  et~al.}{2008}]{Jani2014}
Jani, J.~V., C.~De~Schacht, I.~V. Jani, and G.~Bjune (2008).
\newblock Risk factors for incomplete vaccination and missed opportunity for
  immunization in rural mozambique.
\newblock {\em BMC Public Health\/}~{\em 8}.

\bibitem[\protect\citeauthoryear{Kim, Lim, Cho, and C{\^{o}}t{\'{e}}}{Kim
  et~al.}{2017}]{Kim2017}
Kim, S.~J., G.~J. Lim, J.~Cho, and M.~J. C{\^{o}}t{\'{e}} (2017).
\newblock {Drone-Aided Healthcare Services for Patients with Chronic Diseases
  in Rural Areas}.
\newblock {\em Journal of Intelligent and Robotic Systems: Theory and
  Applications\/}~{\em 88\/}(1), 163--180.

\bibitem[\protect\citeauthoryear{Lee, Assi, Rookkapan, Connor, Rajgopal,
  Sornsrivichai, Brown, Welling, Norman, Chen, Bailey, Wiringa, Wateska, Jana,
  {Van Panhuis}, and Burke}{Lee et~al.}{2011}]{Lee2011}
Lee, B.~Y., T.~M. Assi, K.~Rookkapan, D.~L. Connor, J.~Rajgopal,
  V.~Sornsrivichai, S.~T. Brown, J.~S. Welling, B.~A. Norman, S.~I. Chen, R.~R.
  Bailey, A.~E. Wiringa, A.~R. Wateska, A.~Jana, W.~G. {Van Panhuis}, and D.~S.
  Burke (2011, may).
\newblock {Replacing the measles ten-dose vaccine presentation with the
  single-dose presentation in Thailand}.
\newblock {\em Vaccine\/}~{\em 29\/}(21), 3811--3817.

\bibitem[\protect\citeauthoryear{Lee, Connor, Wateska, Norman, Rajgopal,
  Cakouros, Chen, Claypool, Haidari, Karir, Leonard, Mueller, Paul, Schmitz,
  Welling, Weng, and Brown}{Lee et~al.}{2015}]{Lee2015}
Lee, B.~Y., D.~L. Connor, A.~R. Wateska, B.~A. Norman, J.~Rajgopal, B.~E.
  Cakouros, S.~I. Chen, E.~G. Claypool, L.~A. Haidari, V.~Karir, J.~Leonard,
  L.~E. Mueller, P.~Paul, M.~M. Schmitz, J.~S. Welling, Y.~T. Weng, and S.~T.
  Brown (2015, aug).
\newblock {Landscaping the structures of GAVI country vaccine supply chains and
  testing the effects of radical redesign}.
\newblock {\em Vaccine\/}~{\em 33\/}(36), 4451--4458.

\bibitem[\protect\citeauthoryear{Lee, Haidari, Prosser, Connor, Bechtel,
  Dipuve, Kassim, Khanlawia, and Brown}{Lee et~al.}{2016}]{Lee2016}
Lee, B.~Y., L.~A. Haidari, W.~Prosser, D.~L. Connor, R.~Bechtel, A.~Dipuve,
  H.~Kassim, B.~Khanlawia, and S.~T. Brown (2016, sep).
\newblock {Re-designing the Mozambique vaccine supply chain to improve access
  to vaccines}.
\newblock {\em Vaccine\/}~{\em 34\/}(41), 4998--5004.

\bibitem[\protect\citeauthoryear{Lemmens, Decouttere, Vandaele, and
  Bernuzzi}{Lemmens et~al.}{2016}]{Lemmens2016}
Lemmens, S., C.~Decouttere, N.~Vandaele, and M.~Bernuzzi (2016).
\newblock {A review of integrated supply chain network design models: Key
  issues for vaccine supply chains}.
\newblock {\em Chemical Engineering Research and Design\/}~{\em
  109\/}(February).

\bibitem[\protect\citeauthoryear{Lim, Claypool, Norman, and Rajgopal}{Lim
  et~al.}{2016}]{Lim2016}
Lim, J., E.~Claypool, B.~A. Norman, and J.~Rajgopal (2016, jun).
\newblock {Coverage models to determine outreach vaccination center locations
  in low and middle income countries}.
\newblock ~{\em 9}, 40--48.

\bibitem[\protect\citeauthoryear{Lim, Norman, and Rajgopal}{Lim
  et~al.}{2019}]{Lim2019}
Lim, J., B.~A. Norman, and J.~Rajgopal (2019, dec).
\newblock {Redesign of vaccine distribution networks}.
\newblock {\em International Transactions in Operational Research\/},
  itor.12758.

\bibitem[\protect\citeauthoryear{McCoy and Lee}{McCoy and
  Lee}{2014}]{McCoy2014}
McCoy, J.~H. and H.~L. Lee (2014).
\newblock {Using fairness models to improve equity in health delivery fleet
  management}.
\newblock {\em Production and Operations Management\/}~{\em 23\/}(6), 965--977.

\bibitem[\protect\citeauthoryear{Metcalf, Tatem, Bjornstad, Lessler, O'Reilly,
  Takahashi, Cutts, and GRENFELL}{Metcalf et~al.}{2014}]{Metcalf2014}
Metcalf, C.~J., A.~Tatem, O.~Bjornstad, J.~Lessler, K.~O'Reilly, S.~Takahashi,
  F.~Cutts, and B.~GRENFELL (2014, 08).
\newblock Transport networks and inequities in vaccination: Remoteness shapes
  measles vaccine coverage and prospects for elimination across africa.
\newblock {\em Epidemiology and infection\/}~{\em 143}, 1--10.

\bibitem[\protect\citeauthoryear{{Ministry of Public Health Niger}}{{Ministry
  of Public Health Niger}}{2016}]{cMYP2016}
{Ministry of Public Health Niger} (2016).
\newblock {Comprehensive Multi-Year Plan 2016-2020 on Immunization}.

\bibitem[\protect\citeauthoryear{{National Institute of Statistics of
  Niger}}{{National Institute of Statistics of Niger}}{2012}]{census}
{National Institute of Statistics of Niger} (2012).
\newblock {Census 2012}.
\newblock Available at http://www.citypopulation.de/en/niger, last accessed
  02.14.2021.

\bibitem[\protect\citeauthoryear{Organization}{Organization}{2000}]{SOS}
Organization, W.~H. (2000).
\newblock Sustainable outreach services (sos) : a strategy for reaching the
  unreached with immunization and other services.

\bibitem[\protect\citeauthoryear{Rahmaniani, Crainic, Gendreau, and
  Rei}{Rahmaniani et~al.}{2017}]{rahmaniani2017benders}
Rahmaniani, R., T.~G. Crainic, M.~Gendreau, and W.~Rei (2017).
\newblock The benders decomposition algorithm: A literature review.
\newblock {\em European Journal of Operational Research\/}~{\em 259\/}(3),
  801--817.

\bibitem[\protect\citeauthoryear{Roberts}{Roberts}{2020}]{Berkley2020}
Roberts, L. (2020).
\newblock {Polio, measles, other diseases set to surge as COVID-19 forces
  suspension of vaccination campaigns}.
\newblock
  \url{https://www.science.org/content/article/polio-measles-other-diseases-set-surge-covid-19-forces-suspension-vaccination-campaigns}.
\newblock Accessed 07.19.2021.

\bibitem[\protect\citeauthoryear{Scott and Scott}{Scott and
  Scott}{2017}]{Scott2017}
Scott, J. and C.~Scott (2017, 01).
\newblock Drone delivery models for healthcare.
\newblock In {\em Proceedings of the 50th Hawaii international conference on
  system sciences}.

\bibitem[\protect\citeauthoryear{Toikilik, Tuges, Lagani, Wafiware, Posanai,
  Coghlan, Morgan, Sweeney, Miller, Abramov, Stewart, and Clements}{Toikilik
  et~al.}{2010}]{Toik2010}
Toikilik, S., G.~Tuges, J.~Lagani, E.~Wafiware, E.~Posanai, B.~Coghlan,
  C.~Morgan, R.~Sweeney, N.~Miller, A.~Abramov, A.~Stewart, and C.~J. Clements
  (2010, June).
\newblock Are hard-to-reach populations being reached with immunization
  services? findings from the 2005 papua new guinea national immunization
  coverage survey.
\newblock ~{\em 28\/}(29), 4673—4679.

\bibitem[\protect\citeauthoryear{{UN Office for the Coordination of
  Humanitarian Affairs}}{{UN Office for the Coordination of Humanitarian
  Affairs}}{2018}]{OCHA2018}
{UN Office for the Coordination of Humanitarian Affairs} (2018).
\newblock {Niger - Settlements}.
\newblock Available at https://data.humdata.org/dataset/niger-settlements, last
  accessed 05.21.2020.

\bibitem[\protect\citeauthoryear{UNICEF}{UNICEF}{2018}]{UNICEF2018}
UNICEF (2018).
\newblock {Child given world's first drone-delivered vaccine in Vanuatu}.
\newblock Available at
  https://www.unicef.org/press-releases/child-given-worlds-first-drone-delivered-vaccine-vanuatu-unicef,
  last accessed 02.18.2020.

\bibitem[\protect\citeauthoryear{Vandelaer, Bilous, and Nshimirimana}{Vandelaer
  et~al.}{2008}]{RED}
Vandelaer, J., J.~Bilous, and D.~Nshimirimana (2008).
\newblock {Reaching Every District (RED) approach: a way to improve
  immunization performance}.
\newblock {\em Bulletin of the World Health Organization\/}~{\em 86\/}(3),
  161--240.

\bibitem[\protect\citeauthoryear{Walia, Somarathna, and Jackson}{Walia
  et~al.}{2018}]{Walia2018}
Walia, S., U.~Somarathna, and A.~Jackson (2018).
\newblock {Optimizing the Emergency Delivery of Medical Supplies with Unmanned
  Aircraft Vehicles}.
\newblock In {\em IISE Annual Conference and Expo 2018}, pp.\  1588--1593.

\bibitem[\protect\citeauthoryear{Wedlock, Mitgang, Haidari, Prosser, Brown,
  Krudwig, Siegmund, DePasse, Bakal, Leonard, Welling, Steinglass, Mwansa,
  Phiri, and Lee}{Wedlock et~al.}{2019}]{Wedlock2019}
Wedlock, P.~T., E.~A. Mitgang, L.~A. Haidari, W.~Prosser, S.~T. Brown,
  K.~Krudwig, S.~S. Siegmund, J.~V. DePasse, J.~Bakal, J.~Leonard, J.~Welling,
  R.~Steinglass, F.~D. Mwansa, G.~Phiri, and B.~Y. Lee (2019).
\newblock {The value of tailoring vial sizes to populations and locations}.
\newblock ~{\em 37\/}(4).

\bibitem[\protect\citeauthoryear{WHO}{WHO}{2022}]{WHO2022_jul15}
WHO (2022).
\newblock Covid-19 pandemic fuels largest continued backslide in vaccinations
  in three decades.

\bibitem[\protect\citeauthoryear{{Wingcopter}}{{Wingcopter}}{2021}]{Wingcopter2021}
{Wingcopter} (2021).
\newblock {Technical Details Wingcopter 178 Heavy Lift A}.
\newblock Available at
  https://wingcopter.com/wp-content/uploads/2021/02/Technical-Details-Wingcopter-178-Heavy-Lift-A-Delivery-Variant-1-1.pdf,
  last accessed 08.18.2020.

\bibitem[\protect\citeauthoryear{{World Bank}}{{World
  Bank}}{2020a}]{WorldBank2020}
{World Bank} (2020a).
\newblock {Annual Population Growth Niger}.

\bibitem[\protect\citeauthoryear{{World Bank}}{{World Bank}}{2020b}]{NigerPop}
{World Bank} (2020b).
\newblock {Population Size Niger}.

\bibitem[\protect\citeauthoryear{{World Health Organization}}{{World Health
  Organization}}{2018}]{WHO2018}
{World Health Organization} (2018).
\newblock {Assessment Report of the Global Vaccine Action Plan: Strategic
  Advisory Group of Experts on Immunization}.
\newblock Available at https://apps.who.int
  /iris/bitstream/handle/10665/276967/WHO-IVB-18.11-ara.pdf, last accessed
  05.22.2020.

\bibitem[\protect\citeauthoryear{{World Health Organization}}{{World Health
  Organization}}{2019a}]{WHO2019}
{World Health Organization} (2019a).
\newblock {Children: reducing mortality}.
\newblock Available at https://www.
  who.int/news-room/fact-sheets/detail/children-reducing-mortality, last
  accessed at \\02.17.2020.

\bibitem[\protect\citeauthoryear{{World Health Organization}}{{World Health
  Organization}}{2019b}]{WHO2019b}
{World Health Organization} (2019b).
\newblock {Progress and Challenges with Achieving Universal Immunization
  Coverage}.
\newblock Available at https://www.who.int/immunization/monitoring{\_}
  surveillance/who-immuniz.pdf?ua=1, last accessed 02.21.2020.

\bibitem[\protect\citeauthoryear{{World Health Organization}}{{World Health
  Organization}}{2019c}]{WastageCalc}
{World Health Organization} (2019c).
\newblock {WHO Vaccine Wastage Rates Calculator}.

\bibitem[\protect\citeauthoryear{{World Health Organization}}{{World Health
  Organization}}{2020a}]{WHO2020b}
{World Health Organization} (2020a).
\newblock {Hard fought gains in immunization coverage at risk without critical
  health services, warns WHO}.
\newblock Available at
  https://www.who.int/news-room/detail/23-04-2020-hard-fought-gains-in-immunization-coverage-at-risk-without-critical-health-services-\\warns-who,
  last accessed 05.12.2020.

\bibitem[\protect\citeauthoryear{{World Health Organization}}{{World Health
  Organization}}{2020b}]{Vaccines}
{World Health Organization} (2020b).
\newblock {Immunization Schedule Niger}.

\bibitem[\protect\citeauthoryear{{World Health Organization}}{{World Health
  Organization}}{2020c}]{Vaccines2}
{World Health Organization} (2020c).
\newblock {Prequalified Vaccines}.

\bibitem[\protect\citeauthoryear{Yang, Bidkhori, and Rajgopal}{Yang
  et~al.}{2020}]{Yang2020}
Yang, Y., H.~Bidkhori, and J.~Rajgopal (2020, jan).
\newblock {Optimizing vaccine distribution networks in low and middle-income
  countries}.
\newblock {\em Omega (United Kingdom)\/}.

\bibitem[\protect\citeauthoryear{Zipline}{Zipline}{2019}]{Zipline2019}
Zipline (2019).
\newblock {Zipline,Lifesaving Deliveries by Drone}.
\newblock Available at https://flyzipline.com/, last accessed 02.18.2020.

\end{thebibliography}

\newpage
\section{Appendix}
\subsection{Summary of Notation} \label{AppendixNotation}
%Sets
\begin{table}[htp]
	\centering
	\renewcommand{\arraystretch}{1.5}
	\caption{Sets} 	
	\label{tableSets}
	\begin{tabular}{p{0.15\textwidth}p{0.8\textwidth}}\hline
        \textbf{Set} & \textbf{Description} \\ \hline
		$\mathcal{A}$				&   Arcs\\ 
		$\mathcal{A}^\ast$			&   Arcs that represent land transportation routes, $\mathcal{A}^\ast \subset \mathcal{A}$\\ 
		$\mathcal{A}^\prime$		&   Arcs that represent  community's access to a vaccination center, $\dot{\mathcal{A}}\subset \mathcal{A}$	\\ 
		$\dot{\mathcal{A}}$			&   Arcs that represent potential drone delivery routes,  $\dot{\mathcal{A}}\subset \mathcal{A}$ 	\\ 	
		$\mathcal{I}$				&   Facilities with available cold storage. These are potential locations for drone hubs. \\
		$\mathcal{J}$				&   Cold chain facilities and communities\\
		$\mathcal{K}$				&   Communities, $\mathcal{K}\subset \mathcal{J}$ 	\\ 
		$\bar{\mathcal{K}}$				&   Communities that have access to a clinic,  $\bar{\mathcal{K}} \subset \mathcal{K}$ 	\\ 
		$\mathcal{L}$				&   Vaccines \\ 
		$\mathcal{O}$				&   Outreach posts, $\mathcal{O}\subset\mathcal{J}$\\ 
		$\Pi$ & Vaccination centers, $\Pi\subset\mathcal{J}$\\ 
		$\mathcal{T}$				&   Discrete time periods	\\ \hline		\end{tabular}
    %\end{longtable}
	\renewcommand{\arraystretch}{1}
\end{table}

%Parameter
\begin{table}[htp]
	\centering
	\renewcommand{\arraystretch}{1.5}
	\caption{Parameters} 	\label{tableParameter}
	\begin{tabular}{p{0.15\textwidth}p{0.8\textwidth}}\hline
    \textbf{Parameter} &  \textbf{Description}\\ \hline 
		$a_{l}$						&   Number of doses needed from vaccine $l\in\mathcal{L}$ to complete the vaccination regimen.	\\
		$a_{jk}$						&  An indicator parameter that takes the value 1 if vaccination center $j\in\Pi$ is reachable (within
5Km) from community $k\in\mathcal{K}$, and takes the value 0 otherwise.	\\
		$b$							&   Budget available to establish drone hubs and purchase drones.\\
		$c_{i}$						&   Costs to establish a drone hub at location $i \in {\mathcal{I}}$.\\ 
		$\bar{c}$					&   Unit purchase cost of a drone.\\
		$d_{ij}$					&   Distance (in Km) of arc $(i,j) \in \mathcal{A}$.\\ 
		$\epsilon$ & A scaling factor.\\ 
		$h$							&   Total number of operating hours per drone hub per time period.\\ 
		$q_{l}$						&   Effective packed volume of one dose of vaccine $l\in\mathcal{L}.$   \\
		$m$                         & The maximum number of drones that can be purchased using the available budget $b.$\\ 
		$\pi_{klt}$ & Demand for (number of children in need of a) vaccine $l\in\mathcal{L}$ at community $k\in\mathcal{K}$ in period $t\in\mathcal{T}.$\\
		$\pi_{ilt}$ & Demand for (number of children in need of a) vaccine $l\in\mathcal{L}$ at vaccination center $i\in\Pi$ in period $t\in\mathcal{T}.$\\
		$r_{l}$						&   Volume of the diluent for a single dose of vaccine $l\in\mathcal{L}.$\\
		$s$							&   Average travel speed of a drone.\\
		${u}^d$					&   Drone payload.\\
		${u}_{i}$				&   Cold storage capacity at facility $i \in \mathcal{I}$.\\ 
		${u}^t_{ij}$			&   Cold land transportation capacity along $(i,j)\in\mathcal{A}^{\ast}$.\\ 
		$w^b_{il}$ & Fraction of vaccine $l\in\mathcal{L}$ lost to breakage or inadequate cooling at facility $i\in\mathcal{I}$.\\
        $w^o_{il}$ & Fraction of open-vial wastage during administration of vaccine  $l\in\mathcal{L}$ at vaccination center $i\in\Pi.$ \\
		$w^t_{ijl}$ & Fraction of vaccine $l\in\mathcal{L}$ lost to breakage or inadequate cooling during land transportation along $(i,j)\in\mathcal{A}^{\ast}$.\\
		$w^d_{ijl}$ & Fraction of vaccine $l\in\mathcal{L}$ lost to breakage or inadequate cooling during drone transportation along $(i,j)\in\dot{\mathcal{A}}$.\\
	 \hline
	\end{tabular}
	%\end{longtable}
	\renewcommand{\arraystretch}{1}
\end{table}

%Decision Variables
\begin{table}[htp]
	\centering
	\renewcommand{\arraystretch}{1.5}
	\caption{Decision variables}
	\label{tableDecisionVariables}
	\begin{tabular}{p{0.15\textwidth}p{0.8\textwidth}}\hline
        \textbf{Variable} & \textbf{Description}\\ \hline
		$D_{ijlt}$			&   A continuous variable that represents the number of  doses of vaccine  $l\in\mathcal{L}$ shipped via drones from location $i$ to $j$ ($(i,j)\in\dot{\mathcal{A}}$) in period $t\in\mathcal{T}$.\\
		$I_{ilt}$				&   A continuous variable that represents the inventory of doses of vaccine $l\in\mathcal{L}$ at location $i\in\mathcal{I}$ at the end of time period $t\in\mathcal{T}$.\\
		$N_{k}$						&  A continuous variable that represents the number of fully immunized children of community $k \in \mathcal{K}.$\\
		$S_{ijlt}$			&  A continuous variable that represents the number of doses of vaccine $l\in\mathcal{L}$ shipped via land transportation from facility $i$ to $j$ ($(i,j)\in\mathcal{A}^\ast$) in period  $t\in\mathcal{T}$.\\
		$V_{it}$							&  An integer variable that represents the number of drones used at facility $i\in\mathcal{I}$ during period $t\in\mathcal{T}$.\\
		$W_{i}$	&						A binary variable that takes the value 1 if an outreach post  is established at community $i\in\mathcal{K}$, and takes the value 0  otherwise.\\ 
		$X_{iklt}$			&  A continuous variable that represents the number of doses of vaccine $l\in\mathcal{L}$ administered at vaccination center $i\in\Pi$ to children from community $k\mathcal{K}$	in period $t\in\mathcal{T}$.\\ 
		$X_{ilt}$			&  A continuous variable that represents the number of doses of vaccine $l\in\mathcal{L}$ administered at vaccination center $i\in\Pi$	in period $t\in\mathcal{T}$.\\ 
		$Y_{i}$						&  	A binary variable that takes the value 1 if a drone hub is established at facility $i\in\mathcal{I}$, and takes the value 0  otherwise.\\ 
		$Z$							&  An integer variable that represents the number of drones purchased.\\ \hline
	\end{tabular}
	%\end{longtable}
	%\renewcommand{\arraystretch}{1}
\end{table}

\newpage
\subsection{Model ($Q$)} 
The following is the formulation of model ($Q$).

\begin{align}
%objective
z = max~& \sum_{i \in \Pi} N_{i} + \epsilon  \sum_{i \in \Pi} \sum_{l\in \mathcal{L}}\sum_{t \in \mathcal{T}} X_{ilt}.  \label{FSObj}\\
%budget constraint
s.t.~&\eqref{budget} -\eqref{hubCapacity}, \eqref{initialInventory} -\eqref{storageCapacity}, \eqref{landTransportationCapacity} -\eqref{droneCapacity}, \eqref{binary}-\eqref{Vit} \notag
\end{align}
\begin{eqnarray}
I_{ilt} = (1-w_{il}^b)I_{ilt-1} - \sum_{j:(i,j)\in \mathcal{A}^\ast} S_{ijlt}
 + \sum_{j:(j,i)\in \mathcal{A}^\ast}(1-w_{jil}^t){S}_{jilt-1} - \sum_{j:(i,j)\in \dot{\mathcal{A}}} D_{ijlt} +\notag \\
+ \sum_{j:(j,i)\in \dot{\mathcal{A}}}(1 - {w}^d_{jil})D_{jilt}    - \left(\frac{X_{ilt}}{1-w^o_{il}}\right), \hspace{0.2in} \forall ~i\in\mathcal{I},~l\in \mathcal{L},~t\in \mathcal{T},
\end{eqnarray}
 \begin{eqnarray}
 \sum_{j:(j,i)\in \mathcal{A}^\ast}(1-w_{jil}^t){S}_{jilt-1} 
+ \sum_{j:(j,i)\in \dot{\mathcal{A}}}(1 - {w}^d_{jil})D_{jilt} = \left(\frac{X_{ilt}}{1-w^o_{il}}\right), ~\forall ~i\in\mathcal{O},~l\in \mathcal{L},~t\in \mathcal{T}, \label{inventoryBalance2b} 
\end{eqnarray}
\begin{align}
%demand constraints
&\sum_{t \in \mathcal{T}}{X}_{ilt} \leq \sum_{t \in \mathcal{T}}\pi_{ilt},~&&\forall ~i\in \Pi,~l\in \mathcal{L},\label{demandb} \\
&N_{i} \leq \sum_{t \in \mathcal{T}}\frac{X_{ilt}}{a_{l}},~&&~\forall i \in \Pi, l \in \mathcal{L}.\label{FICb}  
\end{align}

\subsection{Data} 
 
\begin{table}[!h]
\center
\caption{Geographical characteristics of the regions}
\label{tableGeographics}
\begin{tabular}{cccc}\hline
\bf Region & \bf Population {[}-{]} & \bf Area [${km^2}$] & \bf Population density [${km^{-2}}$] \\\hline
Agadez & 626,100   & 634,209 & 0.99                                          \\
Diffa  & 762,700   & 140,216 & 5.44                                          \\
Maradi & 4,728,200 & 38,581  & 122.55                                        \\
Zinder & 4,873,900 & 145,430 & 33.51  \\\hline        
\end{tabular}
\end{table}

\end{document}